\documentclass[10pt]{amsart}
\allowdisplaybreaks
\usepackage[pagewise]{lineno}
\usepackage{fullpage}
\usepackage{amsfonts}
\usepackage{amsmath}
\usepackage{amssymb}
\usepackage{amsthm}
\usepackage{pict2e}
\usepackage{hyperref}
\usepackage{amsbsy}
\usepackage[table]{xcolor}
\usepackage{algorithm}
\usepackage{algpseudocode}
\usepackage{graphicx}
\usepackage{enumerate}
\usepackage{multicol}
\usepackage{multirow}
\usepackage[all,cmtip]{xy}
\usepackage[color=blue!20,textsize=footnotesize]{todonotes}
\usepackage{booktabs}
\usepackage{mathrsfs}
\usepackage{adjustbox}

\newcounter{dummy} \numberwithin{dummy}{section}
\newtheorem{theorem}[dummy]{Theorem}

\newtheorem{lemma}[dummy]{Lemma}

\theoremstyle{remark}
\newtheorem{remark}[dummy]{Remark}
\newtheorem{example}[dummy]{Example}

\newcommand{\scrL}{\mathscr{L}}

\DeclareMathOperator*{\argmax}{argmax}
\DeclareMathOperator*{\argmin}{argmin}

\newcommand{\ve}{\varepsilon}

\DeclareMathOperator{\Ort}{O}

\newcommand{\LibM}{\scrL_{\mathrm{M}}}
\newcommand{\LibMM}{\scrL_{\mathrm{MM}}}
\newcommand{\LibPj}{\scrL_{\mathrm{PolyJet}}}
\newcommand{\LibFIND}{\scrL_{\mathrm{FIND}}}
\newcommand{\LibWS}{\scrL_{\mathrm{WSINDy}}}

\newcommand{\lowbullet}{\raisebox{-0.7ex}{\(\bullet\)}}

\newcommand{\mgpath}{
\begin{adjustbox}{trim=0.1cm 0.1cm 0.1cm 0.1cm}
\begin{tikzpicture}[scale=0.15]
  \node[circle, draw, fill=black, inner sep=1.3pt, outer sep=0pt] (A) at (0,0) {}; %
  \node[circle, draw, fill=black, inner sep=1.3pt, outer sep=0pt] (B) at (2,0) {}; %
  \draw (A) to (B); 
\end{tikzpicture}
\end{adjustbox}
}

\newcommand{\mgloop}{
\begin{adjustbox}{trim=0.1cm 0.1cm 0.1cm 0.1cm}
\begin{tikzpicture}[scale=0.5]
  \node[circle, draw, fill=black, inner sep=1.3pt] (A) at (0,0) {}; %
  \draw (A) to[out=60,in=120,loop] (); %
\end{tikzpicture}
\end{adjustbox}
}

\newcommand{\mglooppath}{
\begin{adjustbox}{trim=0.1cm 0.1cm 0.1cm 0.1cm}
\begin{tikzpicture}[scale=0.5]
  \node[circle, draw, fill=black, inner sep=1.3pt] (A) at (0,0) {}; 
  \node[circle, draw, fill=black, inner sep=1.3pt] (B) at (0.6,0) {}; 
  \draw (A) -- (B);
  \draw (A) to[out=60,in=120,loop] ();
\end{tikzpicture}
\end{adjustbox}
}

\newcommand{\mgpathtwo}{
\begin{adjustbox}{trim=0.1cm 0.1cm 0.1cm 0.1cm}
\begin{tikzpicture}[scale=0.5]
  \node[circle, draw, fill=black, inner sep=1.3pt] (A) at (0,0) {};
  \node[circle, draw, fill=black, inner sep=1.3pt] (B) at (0.6,0) {};
  \node[circle, draw, fill=black, inner sep=1.3pt] (C) at (1.2,0) {};
  \draw (A) to (B);
  \draw (B) to (C);
\end{tikzpicture}
\end{adjustbox}
}

\newcommand{\mgdoublepath}{
\begin{adjustbox}{trim=0.1cm 0.1cm 0.1cm 0.1cm}
\begin{tikzpicture}[scale=0.5]
  \node[circle, draw, fill=black, inner sep=1.3pt] (A) at (0,0) {};
  \node[circle, draw, fill=black, inner sep=1.3pt] (B) at (1,0) {};
  \draw (A) to[bend left=12] (B);
  \draw (A) to[bend right=12] (B);
\end{tikzpicture}
\end{adjustbox}
}

\newcommand{\mgdoubleloop}{
\begin{adjustbox}{trim=0.2cm 0.1cm 0.2cm 0.1cm}
\begin{tikzpicture}[scale=0.5]
  \node[circle, draw, fill=black, inner sep=1.3pt] (A) at (0,0) {};

  \draw (A) to[out=20,in=80,loop] ();
  \draw (A) to[out=100,in=160,loop] ();
\end{tikzpicture}
\end{adjustbox}
}

\numberwithin{equation}{section}

\title[Discovering equivariant PDEs]{Discovering PDEs equivariant under rigid motions}
\author[F.~Ballerin and E.~Grong]{Francesco Ballerin and Erlend Grong}
\date{}

\address{University of Bergen, Department of Mathematics, P.O.~Box 7803, 5020 Bergen, Norway}
\email{francesco.ballerin@uib.no}
\email{erlend.grong@uib.no}

\thanks{The first author is supported by the grant GeoProCo from the Trond Mohn Foundation - Grant TMS2021STG02 (GeoProCo).}

\keywords{Data-driven discovery of PDEs, Equivariance, Sparse regression}

\begin{document}

\begin{abstract}
We consider the problem of PDE discovery from possibly noisy observations under the hypothesis that the underlying dynamic is symmetric in all rigid motions. Rather than using a generic library of derivative monomials, we leverage this assumption to construct libraries whose candidate terms are themselves rigid-motion-equivariant, and combine them with sparse regression to benchmark such libraries over five different equations. The advantage is most pronounced when the noise itself breaks rigid-motion symmetry (e.g., radially or axially varying noise), and when the ambient spatial dimension increases, in which case they are also less resource intensive.
\end{abstract}

\maketitle

\section{Introduction}

PDE discovery is the emerging field of deducing human-readable governing equations from spatial or spatiotemporal data, which usually have some noise. The restriction of human readability usually means that we need a list of discoverable candidate terms, which come in the form of generally nonlinear partial differential operators. In choosing such a library of terms, we would want it to be sufficiently expressive to discover all the terms central to the underlying dynamics. On the other hand, having a larger library of terms risks false positives and will also lead to a higher computational cost. However, it is difficult to make an ``a priori'' argument for which type of terms that are sufficient for discovery. Hence, restrictions of earlier libraries such as in PDE-FIND~\cite{rudy2017data} or in the weak SINDy formulation for PDEs \cite{messenger2021weak} are often imposed based on practical computational concerns. For example, the library of weak SINDy consists of nonlinear polynomial differential operators for which one can take advantage of the integration-by-parts trick when considering the weak formulation of the setting, which is advantageous for computations, especially with noisy observations, but it is hard to see how one would argue from a given physical system that only such terms appear in the unknown dynamics.

One feasible argument to make about a system with unknown dynamics, is to argue that it has certain \emph{symmetries}. In other words, we want to make the a priori assumption that for some self-maps of our space $T: \mathbb{R}^d \to \mathbb{R}^d$, if $u(t,x)$ is a solution of the underlying PDE, then so is $u(t,T(x))$. This assumption is already implicitly present whenever $T$ is a translation in many earlier PDE discovery works, see, e.g., \cite{rudy2017data,messenger2021weak,both2021deepmod,stephany2022pde}, as the candidate terms are assumed to be independent of the spatial variable $x$. In this paper we want to look at PDE discovery for cases where we can assume that our dynamics have \emph{all rigid motions} as symmetries, that is, all maps generated by translations, rotations and reflections. We will consider equations on the form $u_t = N(u)$, where $N$ is a spatial nonlinear partial differential operator satisfying
$$N(u \circ T) = N(u) \circ T, \qquad \text{for any rigid motion $T$.}$$
Such an operator $N$ is called \emph{equivariant under rigid motion} or $E(d)$-\emph{equivariant} where $E(d)$ denotes the Euclidean group of a $d$-dimensional space.
Intuitively, this assumption says that the dynamics does not favor any coordinate system when it comes to center, orientation or order of the axes.

An explicit description of $E(d)$-equivariant nonlinear differential operators for $d \geq 2$ has recently been found in \cite{ballerin2026equivariant}. We want to investigate libraries with such equivariant candidate terms under the following questions.
\begin{enumerate}[$\bullet$]
\item \emph{Recoverability:} For data coming from a PDE that is equivariant under spatial rigid motion, will libraries of equivariant operators more easily recover the correct form of the equation?
\item \emph{Robustness:} While the underlying dynamics may be symmetric under rigid motions, the noise may have a profile that does not follow these symmetries? Will equivariant libraries give more robustness relative to spatial-dependent noise?
\item \emph{Dimension scaling:} It is shown in \cite{ballerin2026equivariant} that the number of equivariant candidate terms with a restricted polynomial degree and maximal number of derivatives in each term is asymptotically constant in dimension, and hence do not grow much with increased dimensions, unlike other libraries that have appeared earlier in literature. Will such libraries perform better in higher dimensions?
\end{enumerate}
To test these questions, we have completed tests on five equations defined by differential operators equivariant under rigid motion: the Allen--Cahn equation, the Swift--Hohenberg equation, the porous medium equation, the Cahn--Hilliard equation and the Kuramoto--Sivashinsky equation. We have looked at the last two in dimension $d=2$, while we for the three first, we have looked at data from dimensions $2\leq d \leq 4$.

Our method of discovery involves numerically computing derivatives from the possibly noisy data and using sparse regression for recovery, similar to, e.g., \cite{brunton2016discovering,rudy2017data,messenger2021weak}. Later works such as \cite{both2021deepmod,stephany2022pde} have shown that including neural networks that approximates both the underlying function $u$ and operator $N$ as an intermediate step can lead to the ability to recover the form of the equation from even more noisy data, but as the final step for these methods also involves sparse regression, we have chosen to only use this for our discovery, as our main aim is to investigate the equivariant libraries. We will have two different libraries consisting of equivariant partial differential operators, both spanning the same space, but having a different basis. We will compare libraries consisting of $E(d)$-equivariant operators with those from PDE-FIND, weak SINDy, as well as the larger library of all polynomials in the jets of a function up to a certain order.

\smallskip

\noindent The main contributions of the paper are the following.
\begin{enumerate}[$\bullet$]
\item We introduce libraries of candidate terms that are equivariant under rigid motions, and leverage them to discover equations equivariant under the same symmetries. Such libraries have sizes that are asymptotically constant in dimension and lead to lower computational cost.. We also show how to modify such libraries to take advantage of integration by parts in a weak formulation of the equations, which is more robust to noise.
\item We demonstrate that our equivariant libraries recover the form of the original equation in more cases in our experiments than comparable non-equivariant alternatives.
\item Our equivariant libraries also do better at recovering the original PDE from the data when considering non-constant noise.
\item We are the first to look at PDE discovery under increasing spatial dimensions. Surprisingly, for all the equations and for all the libraries in our experiments, performance is not affected much by the increase in dimension. However, the equivariant libraries have much lower increase in computational cost.
\end{enumerate}

\subsection{Related work}

The field of PDE discovery has seen intense activity in recent years, with methods spanning linear regression, sparse optimization, weak formulations, symbolic regression, and neural-network-based approaches. All methods are based on the assumption that the solution in the space of coefficients is sparse, i.e., that for an appropriate coordinate system a naturally appearing governing equation is built by a linear combination of few differential operators.

A landmark contribution is the SINDy framework for ordinary differential equations~\cite{brunton2016discovering}, subsequently extended to PDEs as PDE-FIND~\cite{rudy2017data}. For PDE discovery the approach of this class of methods is to construct a library matrix $\mathbf{N}$ of differential operators $N_1, \dots, N_m$ evaluated at a set of collocation points in the spatiotemporal domain. Each column is a differential operator applied to $u$ (in this case monomials $u^n \partial^\alpha u$) with each row being a collocation point. These evaluations can be performed with a variety of numerical methods, such as finite differences or polynomial interpolation. The governing equation $u_t = N(u)$ is then encapsulated in the relationship $\mathbf{U}_t=\mathbf{N}\xi$ with $\mathbf{U}_t$ being the temporal derivative, and the coefficients $\xi_i$, which encode the contribution of the different differential operators. As a minimization problem, we can rewrite this idea as $\xi=\argmin_{w}\lVert \mathbf{U}_t-\mathbf{N} w\rVert_2 + \lambda \lVert w \rVert_0$, with $\lambda\geq 0$ being a cost to push for a small solution, but this minimization problem is unfortunately NP-hard. Therefore the sparse solution is recovered via Sequentially Thresholded Ridge regression (STRidge) instead, a greedy algorithm that iteratively performs ridge regression and prunes small coefficients under the assumption that the true dynamics are sparse in the chosen library. While several canonical equations are correctly recovered, the pointwise formulation is acutely sensitive to noise: numerical differentiation amplifies high-frequency errors, and depending on the governing equation, just a few corrupted evaluations can skew the entire regression.

A significant improvement was introduced by WSINDy~\cite{messenger2021weakODE,messenger2021weak}, which reformulates the regression problem in \emph{weak form}. Firstly, the library is modified to be composed of terms $\partial^\alpha(u^\beta)$, and secondly instead of evaluating the library at individual grid points, each term is integrated against a family of smooth, compactly supported test functions $\psi_{ij}$. The linear nature of the terms allows the derivatives to be transferred onto the test functions via integration by parts:
\begin{equation}
  \int \psi_{ij}\,\partial^\alpha(u^\beta)\,dx
  = (-1)^{|\alpha|} \int \partial^\alpha\psi_{ij}\,u^\beta\,dx, \qquad \text{for any multi-index $\alpha$} .
\end{equation}
This has two beneficial effects: first, high-order derivatives of the (possibly noisy) data $u$ are replaced by analytically known derivatives of the smooth test functions, substantially reducing noise amplification; second, the convolution with $\psi_{ij}$ provides an intrinsic local averaging that further suppresses measurement noise.  WSINDy also introduces the Modified Sequential Thresholding Least Squares algorithm (MSTLS), a variation of sequential thresholding tailored to the weak-form setting.

In parallel, a substantial body of work tackles PDE discovery through deep learning.  Physics-informed neural networks (PINNs) can encode governing equations in the loss function, enabling simultaneous solution and discovery of PDEs \cite{raissi2019physics}, while PDE-Net~\cite{long2018pdenet} learns differential operators as constrained convolution kernels.  Although powerful, these methods generally do not return a human-readable symbolic expression without a separate symbolic post-processing step, placing them in a different category from the library-based sparse regression pursued here. Example of post-processing step is performed by DeepMoD~\cite{both2021deepmod} and PDE-READ \cite{stephany2022pde}. After training a neural network, the latter uses such data to run sparse regression with an algorithm named RFE. This introduces a user-defined importance metric which yields a sequence of progressively sparser candidate solutions by identifying at each step the least important feature to be removed. There are also papers taking equivariance into account when learning or solving PDEs with deep learning. The paper \cite{he2022neural} considers PDEs of vector valued functions that are equivariant under all affine transformations, while \cite{lagrave2022equivariant} gives a general overview of equivariant neural nets for solving PDEs.

\section{Problem and underlying assumptions}
The objective of this paper is to discover a (generally nonlinear) partial differential operator $N: f \mapsto N(f)$ defined on functions $f:\mathbb{R}^d \to \mathbb{R}$, which we consider as the $d$ spatial dimensions. We want to discover this operator from (possibly noisy) spatiotemporal data $\mathbf{U} = \{ U(t_k,x_k)\}_{k=1}^K = \{ u(t_k, x_k) +\ve(t_k,x_k) \}_{k=1}^K$ of a function $u: [0,T] \times \mathbb{R}^d \to \mathbb{R}$ solving the equation $u_t = \partial_t u = N(u)$, with $\ve$ denoting the noise. Our underlying assumption is that we have a library $\scrL = \{ N_1, \dots, N_m\}$ of differential operators such that the operator $N$ that we are interested in can be written as a linear combination of terms in this library;
\begin{equation} \label{Nspan} N = \sum_{j=1}^m \xi_j N_j, \quad \xi = (\xi_1, \dots, \xi_m) \in \mathbb{R}^m.\end{equation}
Using numerical derivatives and our data $\mathbf{U}$, we obtain a vector $\mathbf{U}_t$ and an array $\mathbf{N} = \mathbf{N}(\mathbf{U},\scrL)$ computing respectively $U_t$ and $(N_1(U), \dots, N_m(U))$ at our given sample times and points. We then want to discover $N$ by solving
$$\mathbf{U}_t = \mathbf{N} \xi, \qquad \mathbf{U}_t \in \mathbb{R}^K, \mathbf{N} \in \mathbb{R}^{K \times m}, \xi \in \mathbb{R}^m.$$
We will consider this problem under the following assumptions.

\smallskip
\smallskip

\paragraph{\it Sparsity:} It is a principled assumption for a differential operator $N$ modeling a natural system to be sparse. As outlined in \cite{brunton2016discovering}, nature often tends to prefer simple equations and symmetry. It is not unreasonable that the same applies here, so that $N$ as in \eqref{Nspan} can be written as a linear combination of \emph{only a few of the terms in the library $\scrL = \{ N_1, \dots, N_m\}$}, meaning that the vector $\xi$ is sparse. For the rest of the paper, let
$$\textstyle \|U\|_s = \left(\sum_{j=1}^m |U_j|^s \right)^{1/s}, \qquad \|U\|_0 = \left|\{ j \, : \, U_j \neq 0\}\right|, \qquad s \geq 1, \quad U \in \mathbb{R}^m.$$
An estimate $w= \hat \xi$ to the true coefficients $\xi$ in \eqref{Nspan} can be easily obtained by several methods, among which ordinary least square (OLS) regression by $w = \argmin_{\tilde w} \| \mathbf{U}_t - \mathbf{N}\tilde w\|_2^2$ is the most common. However, the columns in $\mathbf{N}$ will typically be highly correlated (see Section~\ref{sec:Covariance}), leading to high variance in the outcome of the OLS regression and leaving it vulnerable to noise and numerical errors. One can lower this variance at the cost of having some bias by instead completing a ridge regression, see, e.g., \cite[Chapter~6.2]{james2013introduction}. Ridge regression looks for $w = \argmin_{\tilde w} (\| \mathbf{U}_t - \mathbf{N} \tilde w\|_2^2 + \lambda \|\tilde w\|^2_2)$ for some $\lambda \geq 0$. Increasing $\lambda$ will increase bias and decrease variance. Note that for any $\lambda >0$, there is an $r = r(\lambda) >0$ such that the outcome of the ridge regression equals $\argmin_{\|\tilde w\|_2 \leq r} \|\mathbf{U}_t - \mathbf{N} \tilde w \|_2^2$. However, neither OLS, nor ridge regression give sparse solutions. In order to obtain the desired sparse estimate, we could look at the solution of
\begin{equation} \label{lzeropunish} \argmin_{\tilde w} \left(\| \mathbf{U}_t - \mathbf{N} \tilde w \|_2^2 + \lambda \| \tilde w \|_0 \right)\end{equation}
with $\lambda \geq 0$ being some cost on the lack of sparsity, but such a minimization is NP-hard \cite[Theorem~1]{natarajan1995sparse}. The standard alternative approach is to use LASSO (Least Absolute Shrinkage and
Selection Operator) regression, which means finding $\argmin_{\tilde w} \left(\| \mathbf{U}_t - \mathbf{N} \tilde w \|_2^2 + \lambda \| \tilde w \|_1 \right)$ \cite{tibshirani1996regression}, but this approach also has problems when the columns of $\mathbf{N}$ have high correlations \cite{herawati2018regularized}. We will therefore use algorithms STRidge and RFERidge, see Section~\ref{sec:Regression}.

\smallskip
\smallskip

\paragraph{\it Polynomial differential operator:} We assume that the nonlinear operator $N$ is \emph{polynomial}. In other words, we assume that $N(f)$ can be written as a polynomial in $f$ and its iterated derivatives. The following terminology will also be useful. Consider a partial differential operator $P$ defined by
\begin{equation} \label{monomial} P(f) = (L_1 f)(L_2f) \cdots (L_n f),\end{equation}
where $L_j$ is a linear partial differential operator of order $r_j$ for $j=1, \dots, n$. We call $n$ \emph{the polynomial degree} of $P$ as in \eqref{monomial}, $r = \max_j \{r_j\}$ its \emph{order} and $p = \sum_{j=1}^n r_j$ its \emph{total order}. Any polynomial differential operator~$N$ can be written as a linear combination $\sum_{l=1}^{l_{\max}} P_l$ of operators on the form \eqref{monomial}, with polynomial degree, order and total order defined as the maximal of that of $P_1$, $\dots$, $P_{l_{\max}}$. This is equivalent to stating that~$N(f)$ can be written as a polynomial in degree $n$ in $f$ and its iterated derivatives, such that the highest order derivative of $f$ used is order $r$ and where each term has a maximal of $p$ derivatives when taking all factors into consideration.
\emph{We will assume that $N$ is a polynomial differential operator with an assumed bound for the polynomial degree and total order}.

\smallskip
\smallskip

\paragraph{\it Equivariance with respect to rigid motions}
\emph{A rigid motion} or \emph{isometry} of $\mathbb{R}^d$ is a function $T:\mathbb{R}^d \to \mathbb{R}^d$ such that $\| T(x) - T(y)\|_2 =\|x-y\|_2$ for any $x,y \in \mathbb{R}^d$. All such maps can be written as $T(x) = Ax+b$ where $b \in \mathbb{R}^d$ and $A \in \Ort(d)$ is an orthogonal matrix. All of these maps collectively form the Euclidean group $E(d)$. We will assume that the operator $N$ is \emph{equivariant with respect to rigid motion} in the sense that for any such~$T \in E(d)$ and smooth $f: \mathbb{R}^d \to \mathbb{R}$, we have that
$$N(f \circ T) = N(f) \circ T.$$
We will take advantage of the structure found using the results found in  \cite{ballerin2026equivariant}.

\begin{theorem} \label{th:multigraph}
Let $d \geq 2$ be fixed. There is a surjective map from the vector space spanned by all isomorphism classes of multigraphs with $\leq n$ vertices and $\leq p$ edges on the space of polynomial nonlinear differential operators with polynomial degree $\leq n$ and total order $\leq 2p$ on $\mathbb{R}^{d}$. Furthermore, if $p \leq d$, the map is bijective. 
\end{theorem}
We refer to \cite{ballerin2026equivariant} for the details on multigraphs and the correspondence. We remark that the same reference also shows that nonlinear operators equivariant under rigid motions always have an even total degree. Also, let $\mathscr{V}(n,p,\mathbb{R}^d)$ be the vector space spanned by equivariant polynomial nonlinear differential operators of polynomial degree at most $n$ and total order at most $2p$. Then Theorem~\ref{th:multigraph} says that the dimension of this space has as upper bound the number of non-isomorphic multigraphs with at most $n$ and $p$ edges. Furthermore, this bound is realized whenever $p \leq d$, making $\dim \mathscr{V}(n,p,\mathbb{R}^d)$ asymptotically constant in $d$.

\begin{example}
We will use the following restrictions for the libraries in our paper. We will assume that we are searching for an operator $N$ of polynomial degree $\leq 3$ and total order $\leq 4$. Since $4 = 2 p$ with $p=2$, we can define a one-to-one correspondence between the vector space spanned by multigraphs with $\leq 3$ vertices and $\leq 2$ edges. There are exactly 20 isomorphism classes of such multigraphs, so the number of equivariant partial differential operators of polynomial degree $\leq 3$ and total order $\leq 4$ is 20 in any dimension $d$ greater or equal to $2$. These are all listed in Table~\ref{tab:multigraph_operators}.
\end{example}

\begin{table}[htbp]
  \centering
  \renewcommand{\arraystretch}{1.5} 
  
  \small
  \begin{tabular}{|r|cccc||cccc|}
  \hline
    $p$  & $n=0$ & $n=1$ & $n=2$ & $n=3$ & $n=0$ & $n=1$ & $n=2$ & $n=3$    \\  \hline
    $0$  & $\emptyset$ & $\lowbullet$ & $\lowbullet$ $\lowbullet$ & $\lowbullet$ $\lowbullet$ $\lowbullet$ & $1$ & $f$ & $f^2$ & $f^3$  \\ \hline
    $1$  & & \mgloop & $\lowbullet$ \mgloop & $\lowbullet$ $\lowbullet$ \mgloop & & $\Delta f$ & $f \Delta f$ & $f^2 \Delta f$ \\
    &  & & \mgpath & $\lowbullet$ \mgpath & & & $\|\nabla f\|^2$ & $f \| \nabla f \|^2$  \\\hline
    $2$ & & \mgdoubleloop & $\lowbullet$ \mgdoubleloop & $\lowbullet$ $\lowbullet$ \mgdoubleloop & & $\Delta^2 f$ & $f \Delta^2 f$ & $f^2 \Delta^2 f$ \\
    &  & & \mglooppath & $\lowbullet$ \mglooppath & & & $\langle \nabla \Delta f , \nabla f \rangle$ & $ f\langle \nabla \Delta f , \nabla f \rangle$ \\
    & &  & \mgdoublepath & $\lowbullet$ \mgdoublepath & & & $\| \nabla^2 f \|^2$ & $f\| \nabla^2 f\|^2$ \\
    &  & & \mgloop \mgloop  & $\lowbullet$ \mgloop \mgloop & & & $(\Delta f)^2$ & $f (\Delta f)^2$ \\
    &  & & & \mgpathtwo & & & &  $\nabla^2 f(\nabla f, \nabla f)$   \\
    &  & & & \mgpath \mgloop & & & & $\| \nabla f \|^2  \Delta f$ \\
    \hline
  \end{tabular}
  \caption{The 20 $\Ort(d)$-equivariant polynomial nonlinear differential operators
           on $\mathbb{R}^d$ of total differential order~$2p$ ($p \le 2$) and
           polynomial degree~$n$ ($n \le 3$). These are listed with the corresponding isomorphism classes of multigraphs with up to 3 vertices and up to 2 edges. Notice that loops and multiple edges between vertices are allowed in multigraphs. We will call the library of terms of this table \emph{the multigraph library} and write it as $\LibM$.}
  \label{tab:multigraph_operators}
\end{table}

\section{Methodology}
The objective is to learn a partial differential operator $N$ that is equivariant with respect to rigid motions by estimating equation coefficients $\xi$ as in \eqref{Nspan}, given spatio-temporal data $\mathbf{U}$ possibly affected by noise. In pursuing such objective, we compare libraries that only consist of equivariant operators (as introduced in \cite{ballerin2026equivariant}) to well-established libraries and methods in PDE discovery. The coefficient vector $\xi$ will be estimated under the prior that the solution is sparse in terms of coefficients, and therefore sparse regression is employed. The data can be noisy, and using a weak formulation will help us deal with this noise.

\subsection{The chosen equations}
Our results come from experiments with respect to five equations involving operators that are equivariant under rigid motions.
\begin{enumerate}[1.]
\item The \emph{Allen--Cahn} equation is a scalar reaction-diffusion equation introduced by Allen and Cahn~\cite{allen1979microscopic} to model the motion of antiphase boundaries in crystalline solids, with relations to minimal surfaces \cite{savin2010phase,chodosh2020minimal}. It can be written as 
\begin{equation} \tag{AC} \label{AC}
    u_t = \varepsilon^2 \Delta u + u - u^3, \qquad \ve = 0.15,
\end{equation}
with polynomial degree $3$ and total order $2$. In our benchmark suite it is the simplest equation, containing only three active terms, all of which are pure monomials with no mixed derivative-nonlinear couplings. The parameter $\varepsilon$ controls the width of the diffuse interface separating the two stable phases.
Simulations are initialized from small-amplitude Gaussian noise ($0.05\,\mathcal{N}(0,1)$ per grid point) on the periodic domain $[0,2\pi]^2$ ($64^2$ grid) and evolved with an IMEX pseudo-spectral scheme. Snapshots of a solution are shown in Figure~\ref{fig:ac_motion}.

\begin{figure}
    \centering
    \includegraphics[width=0.85\linewidth]{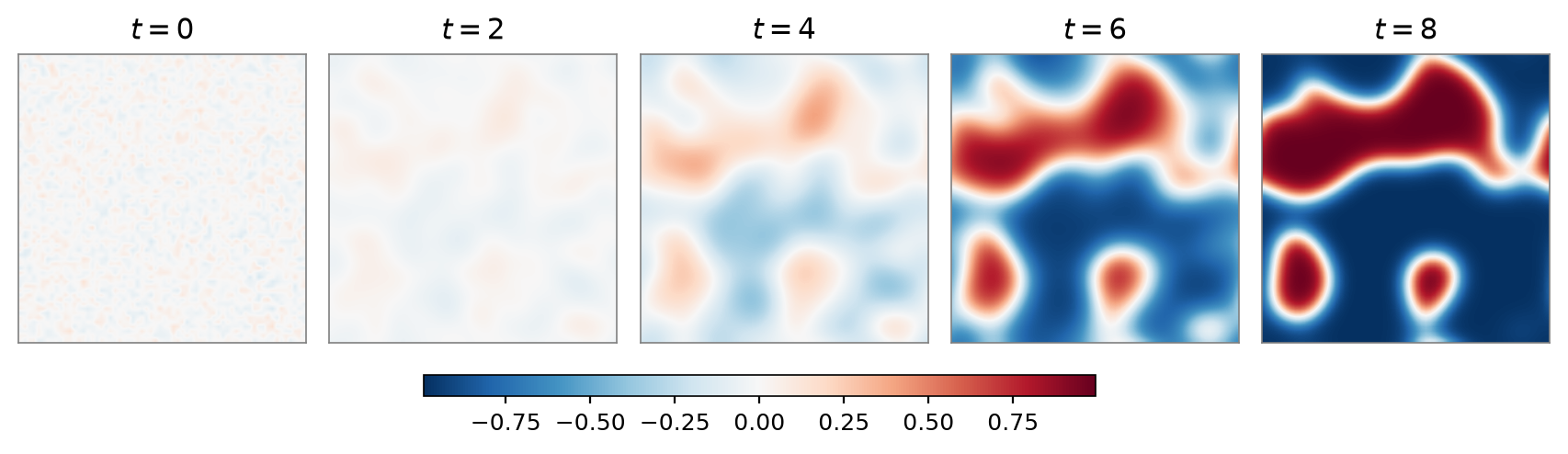}
    \caption{Snapshots of the evolution of the Allen--Cahn equation with $\varepsilon = 0.15$ on the periodic domain $[0,2\pi]^2$ ($64^2$ grid).}
    \label{fig:ac_motion}
\end{figure}

\item \emph{The Swift--Hohenberg equation} was introduced by Swift and Hohenberg~\cite{swift1977hydrodynamic} as a normal form for the onset of convective instability in Rayleigh--Bénard cells. It has since become a model for pattern formation, i.e. the spontaneous emergence of spatially periodic structures (stripes, hexagons, squares) in dissipative systems far from equilibrium~\cite{cross1993pattern}. The equation is
\begin{equation} \tag{SH} \label{SH}
    u_t = (r-1)u - 2\Delta u - \Delta^2 u - u^3, \qquad r = 0.3,
\end{equation}
with $r$ controlling linear stability. For $r>0$ the system self-organizes into steady patterns. Simulations are initialized from small-amplitude Gaussian noise ($0.01\,\mathcal{N}(0,1)$ per grid point) on the periodic domain $[0,8\pi]^2$ ($64^2$ grid) and evolved with an IMEX pseudo-spectral scheme. Snapshots of a solution are shown in Figure~\ref{fig:sh_motion}.

\begin{figure}
    \centering
    \includegraphics[width=0.85\linewidth]{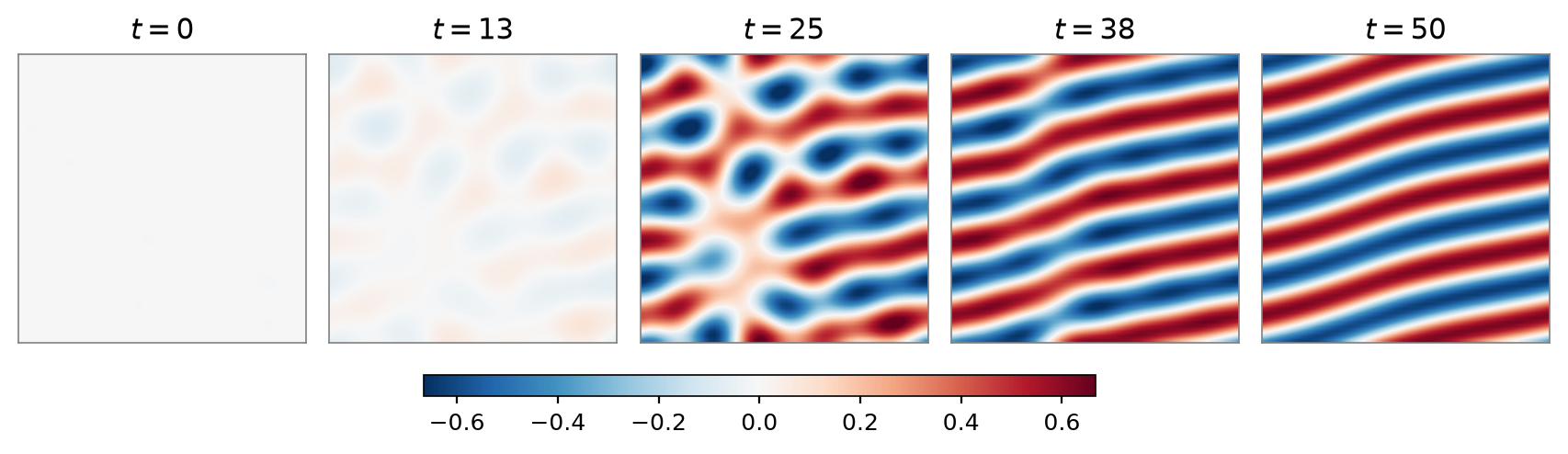}
    \caption{Snapshots of the evolution of the Swift--Hohenberg equation with $r = 0.3$ on $[0,8\pi]^2$ ($64^2$ grid). Steady patterns are visible in the form of stripes, one of possible patterns for this equation.}
    \label{fig:sh_motion}
\end{figure}

\item \emph{The porous medium equation} is the class of nonlinear diffusion equations $u_t = \Delta(u^m)$, which model gas flow through porous media, groundwater infiltration, and other processes where the diffusion coefficient depends on the concentration itself~\cite{vazquez2006porous}. We take the case $m = 3$, which in divergence form reads
\begin{equation} \tag{PM} \label{PM}
    u_t = \frac{1}{3} \Delta u^3 = \nabla \cdot (u^2 \nabla u) = u^2 \Delta u + 2u\|\nabla u\|^2.
\end{equation}
with polynomial degree $3$ and total differential order $2$. The factor $1/3$ is inessential and can be absorbed by rescaling time. Simulations use an explicit Euler scheme with spectral derivatives on $[0,6\pi]^2$ ($64^2$ grid), initialized from a centered Gaussian distribution of width $L/8$ (normalized to unit peak). Snapshots of the spreading are shown in Figure~\ref{fig:pme_motion}.

\begin{figure}
    \centering
    \includegraphics[width=0.85\linewidth]{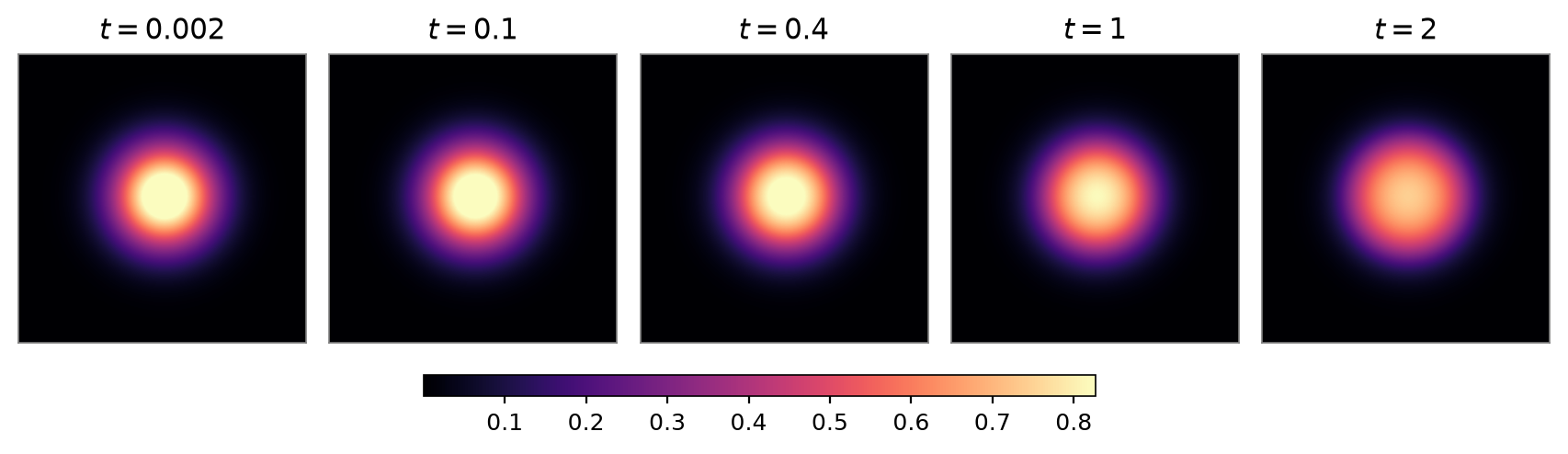}
    \caption{Snapshots of the spreading of a compactly supported initial mass under the porous-medium equation on $[0,6\pi]^2$ ($64^2$ grid).}
    \label{fig:pme_motion}
\end{figure}

\item \emph{The Cahn--Hilliard equation} \cite{cahn1958free} governs phase separation in a binary alloy, and more generally describes spinodal decomposition in systems with a conserved order parameter. Applications also include diblock copolymers, image inpainting, and tumor-growth models~\cite{miranville2019cahn,kim2016basic}. The equation has polynomial degree $3$ and total differential order $4$, written compactly as

\begin{equation} \tag{CH} \label{CH}
    u_t = -\varepsilon^2 \Delta^2 u + \Delta(u^3 - u)
    = -\Delta u + 3u^2 \Delta u + 6u\|\nabla u\|^2 - \varepsilon^2 \Delta^2 u,
    \qquad \varepsilon = 0.2.
\end{equation}
where $\varepsilon$ controls the thickness of the interface between the separating phases. Simulations are initialized from Gaussian noise of amplitude $0.5$ on the periodic domain $[0,8\pi]^2$ ($128^2$ grid) and evolved with an IMEX pseudo-spectral scheme. Snapshots of a solution are shown in Figure~\ref{fig:ch_motion}.

\begin{figure}
    \centering
    \includegraphics[width=0.85\linewidth]{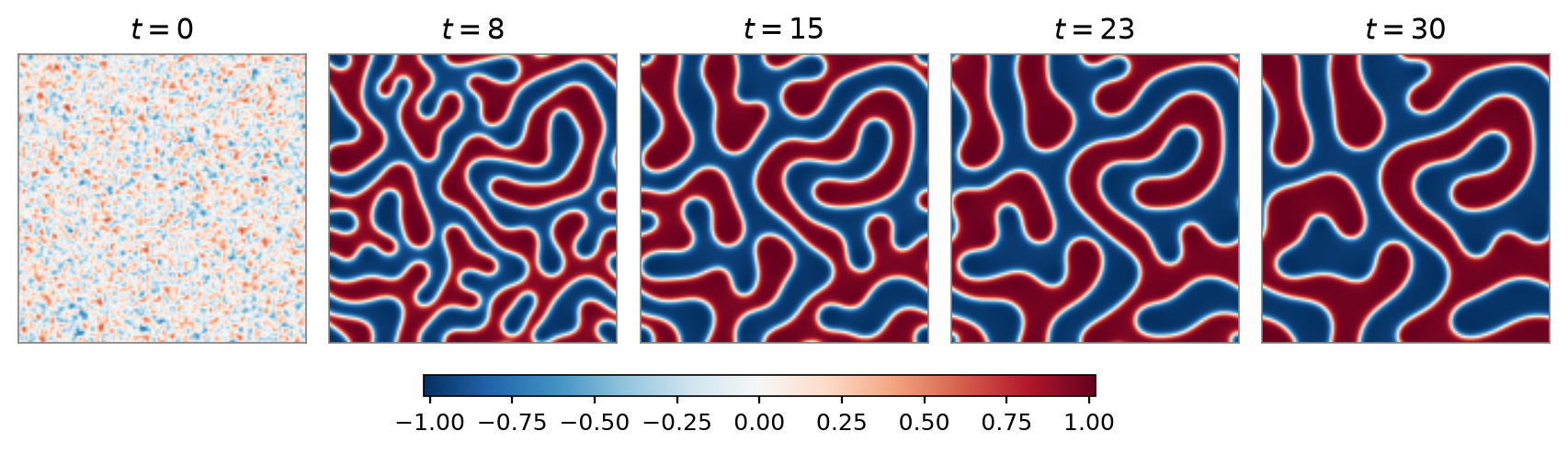}
    \caption{Snapshots of evolution of the Cahn--Hilliard equation with $\varepsilon = 0.2$ on $[0,8\pi]^2$ ($128^2$ grid).}
    \label{fig:ch_motion}
\end{figure}

\item \emph{The Kuramoto--Sivashinsky (KS) equation} was derived independently by Kuramoto and Tsuzuki~\cite{kuramoto1976persistent} in the context of reaction-diffusion systems, and by Michelson and Sivashinsky~\cite{michelson1977nonlinear} in the study of flame-front instabilities, and it is a prototypical model of \emph{spatiotemporal chaos}. The complexity of the dynamics grows with the domain size~\cite{hyman1986kuramoto,michelson1986steady} and can be categorized in \emph{stable}, \emph{cellular}, and \emph{chaotic} regimes. The equation has total differential order $4$ and polynomial degree $2$, given by
\begin{equation} \tag{KS} \label{KS}
    u_t = -\tfrac{1}{2}\|\nabla u\|^2 - \Delta u - \Delta^2 u
\end{equation}
We simulate three regimes: cellular ($L = 4\pi$), weakly chaotic ($L = 8\pi$), and strongly chaotic ($L = 16\pi$). We exclude the stable regime from the analysis as it appears uninteresting, due to its dynamic decaying to zero. All simulations use the ETDRK4 exponential integrator~\cite{cox2002exponential,kassam2005fourth} initialized from Gaussian noise of amplitude $0.1$.  Snapshots of a solution are displayed in Figure~\ref{fig:ks_motion}. 

\begin{figure}
    \centering
    \includegraphics[width=0.85\linewidth]{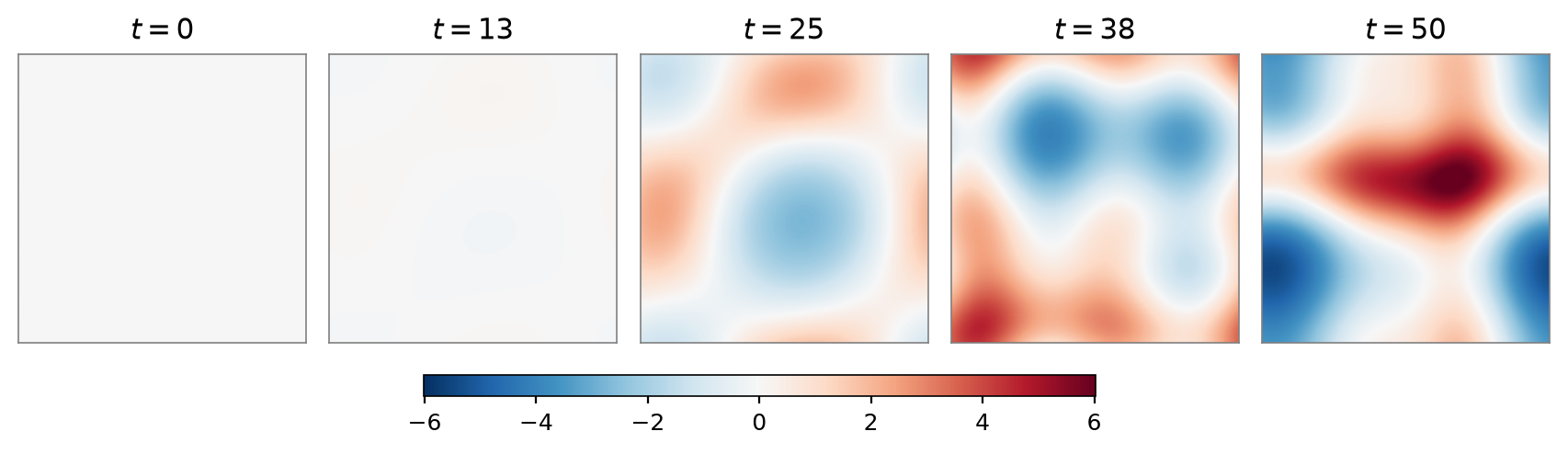}
    \includegraphics[width=0.85\linewidth]{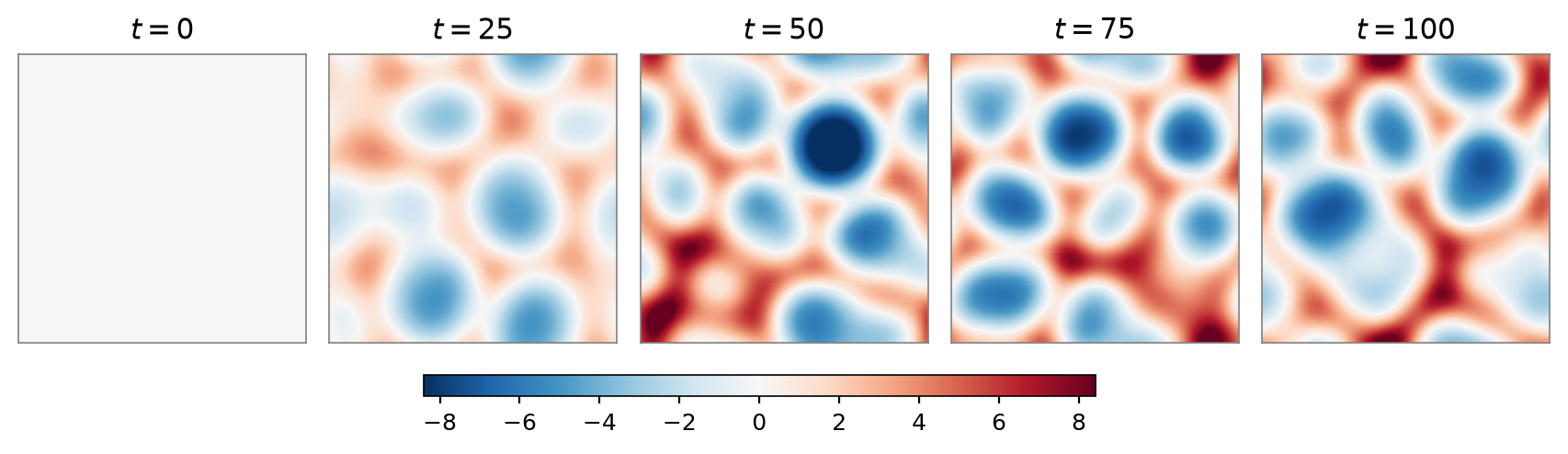}
    \caption{Snapshots of evolution of the Kuramoto--Sivashinsky equation in cellular regime (top) and chaotic regime (weakly chaotic, bottom). Cellular regime is attained on a periodic domain $[0,4\pi]^2$, while chaotic regime is attained on a $[0,8\pi]^2$ domain, both on a $64^2$ grid.}
    \label{fig:ks_motion}
\end{figure}
\end{enumerate}
We will consider all equations in experiments with $d=2$ in Section~\ref{sec:ResultsTwoD} and experiments with different noise profiles in Section~\ref{sec:SpaceVarying}. We will consider the Allen--Cahn equation, Swift--Hohenberg equation, and the porous medium equation in higher dimensions in Section~\ref{sec:IncreasedDimension}.

\subsection{Data and sampling} \label{sec:Data}
All our experiments are run on synthetic data, so that the operator $N$ we try to discover is known exactly from the start. We work on a flat domain with periodic boundary conditions, discretized by a uniform grid of $N_{\mathrm{grid}}^d$ points with spacing $\Delta x = L/N_{\mathrm{grid}}$. Periodicity means that all libraries and both derivative backends can be evaluated at every grid point without any special treatment near a boundary, so no method is advantaged by its boundary handling.

Solutions are computed by pseudo-spectral time-stepping. For the Allen--Cahn, Swift--Hohenberg, and Cahn--Hilliard equations we use an IMEX scheme, treating the linear part implicitly and the nonlinear part explicitly. For the porous media equation we resort to explicit Euler, as it contains no part to be treated implicitly. The Kuramoto--Sivashinsky equation is much stiffer, and is instead integrated with the fourth order exponential time-differencing Runge--Kutta scheme ETDRK4 \cite{cox2002exponential, kassam2005fourth}. Initial conditions are random and drawn from a fixed seed. The integration parameters are collected in Table~\ref{tab:simulation_parameters}. The Kuramoto--Sivashinsky equation appears three times corresponding to the three regimes (cellular, weakly chaotic, strongly chaotic).

\begin{table}[htbp]
  \centering
  \small
  \begin{tabular}{l|cccccl}
    Equation & $N_{\mathrm{grid}}$ & $L$ & $T$ & $\Delta t$ & $n_{\mathrm{save}}$ & Parameters \\ \hline
    Allen--Cahn \eqref{AC}                & $64$  & $2\pi$  & $8$   & $0.01$  & $20$ & $\ve = 0.15$ \\
    Swift--Hohenberg \eqref{SH}           & $64$  & $8\pi$  & $50$  & $0.1$   & $5$  & $r = 0.3$ \\
    Cahn--Hilliard \eqref{CH}             & $128$ & $8\pi$  & $30$  & $0.005$ & $30$ & $\ve = 0.2$ \\
    Porous medium \eqref{PM}              & $64$  & $6\pi$  & $2$   & $0.002$ & $25$ & $m = 3$ \\
    Kuramoto--Sivashinsky \eqref{KS} (cellular)       & $64$  & $4\pi$  & $50$  & $0.1$   & $5$  & \\
    Kuramoto--Sivashinsky \eqref{KS} (weakly chaotic) & $64$  & $8\pi$  & $100$ & $0.1$   & $5$  & \\
    Kuramoto--Sivashinsky \eqref{KS} (strongly chaotic)  & $64$  & $16\pi$ & $100$ & $0.1$   & $5$  & \\
  \end{tabular}
  \caption{Simulation parameters in dimension $d=2$. The solution is computed on $[0,L]^2$ with time step $\Delta t$ up to time $T$, and every $n_{\mathrm{save}}$-th state is retained, so that the spacing between two saved snapshots is $\Delta t_{\mathrm{snap}} = n_{\mathrm{save}} \Delta t$.}
  \label{tab:simulation_parameters}
\end{table}

\subsection{Sparse regression and training of hyperparameters} \label{sec:Regression}
We will apply two main methods for sparse regression. The first approach is STRidge \cite{brunton2016discovering,rudy2017data}. After normalizing columns of $\mathbf{N}$ to have unit length, we complete repeated ridge regressions while removing coefficients that are under a given tolerance level $tol >0$ in magnitude for each iteration until the number of coefficients stabilizes. This approach has two hyperparameters in the tolerance level $tol$ and a parameter $\lambda$ for the ridge regression. We find the correct tolerance level $tol >0$ by splitting our data into training and test data, completing ridge regression on the training data, and adjust the tolerance level according to the solution's performance on the test data relative to a cost as in \eqref{lzeropunish}.

We have also tested a modification to STRidge introduced in \cite{messenger2021weak}, called modified sequential-thresholding least-squares (MSTLS). It also involves iterated regressions, but rather than removing coefficients that fall below a certain tolerance level, MSTLS removes coefficients whose magnitude fall outside a permitted interval in iterated OLS regressions. Relative to a given threshold $\mu >0$, the coefficient $|w_j|$ needs to be above $\mu \cdot \max\{ 1, \|\mathbf{U}_t\|_2/\|\mathbf{N}_j\|_2 \}$ and below $\frac{1}{\mu} \cdot \min\{ 1, \|\mathbf{U}_t\|_2/\|\mathbf{N}_j\|_2 \}$. If $w^\mu$ is the result of this iteration and $w^0$ is the outcome of OLS, then the optimal $\mu$ is chosen according to the cost
\begin{equation} \label{muMin}\mu \mapsto \frac{\| \mathbf{N} (w^\mu-w^0)\|_2}{\|\mathbf{N} w^0\|_2} + \frac{\|w^\mu\|_0}{m}.\end{equation}
Minimizing is done by choosing a grid for the $\mu$-values and choosing the minimal value.

More recently, the method of Recursive Feature Elimination (RFE) \cite{guyon2002gene} has been used for PDE-discovery in \cite{stephany2022pde}. It completes iterative ordinary least square regression, removing the coefficient that has least impact in each iteration. To be more precise, if $e_j$, $j=1, \dots, m$ is the $j$-th unit vector in $\mathbb{R}^m$, then a coefficient $w_j \neq 0$ of $w = (w_1, \dots, w_m)$ is considered to have least impact if $j = \argmin_{J, w_J \neq 0} \|\mathbf{U}_t - \mathbf{N} (w - w_J e_J)\|^2_2$. If $w^{(0)} = w$ is the result of an initial OLS regression and $w^{(j)}$ is the result after completing $j$-th iteration of completing an OLS regression and then setting the least impactful coefficient to zero, we write $Res^{(j)} = \| \mathbf{U}_t - \mathbf{N} w^{(j)}\|_2^2$. By this definition $Res^{(j)}$ is increasing in $j$ for $j=0,1, \dots, m$ with $\xi^{(m)} =0$ and $Res^{(m)} = \|\mathbf{U}_t\|_2^2$. For $j=0, \dots, m-1$, we then define $Ranking^{(j)} = Res^{(j+1)}/Res^{(j)}$. RFE then gives out $w^{(j)}$ if $Ranking^{(j)}$ is maximal. Note that RFE does not require any hyperparameters. We apply a generalization of this approach RFE, which we will call RFERidge, where we complete iterated ridge regressions instead, again to deal with the correlations of the columns. This introduces a hyperparameter $\lambda$, and we now have residual size $Res^{(j)} = \| \mathbf{U}_t - \mathbf{N} w^{(j)} \|_2^2 + \lambda \| w^{(j)}\|^2_2$. We recover the usual RFE regression by setting $\lambda =0$.

\begin{algorithm}[H]
\caption{STRidge$(y,X_0, \lambda, tol)$}
\label{alg:STRidge}
\begin{algorithmic}[1]
\Require $y \in \mathbb{R}^K$, $X_0 \in \mathbb{R}^{K \times m}$, $\lambda \geq 0$, $tol > 0$, $y \neq 0$
\State $Scale_j \gets 1/\|X_0[:,j]\|_2$, \; $X \gets X_0 \cdot Scale$ \Comment{normalize columns}
\State $small \gets \emptyset$, \; $newsmall \gets \{0\}$, \; $w \gets 0 \in \mathbb{R}^m$
\While{$|newsmall| > |small|$}
    \State $small \gets newsmall$, \; $big \gets \{0,\dots,m-1\} \setminus small$
    \State $w[big] \gets \argmin_\psi (\|y - X[:,big]\,\psi\|_2^2 + \lambda\|\psi\|_2^2)$, \quad $w[small] \gets 0$
    \State $newsmall \gets \{i : |w_i| < tol\}$
\EndWhile
\State $w[big] \gets \argmin_\psi \|y - X[:,big]\,\psi\|_2^2$ \Comment{unbias on the support}
\State \Return $w \cdot Scale$
\end{algorithmic}
\end{algorithm}

\begin{algorithm}[H]
\caption{RFERidge$(y,X_0, \lambda)$}
\label{alg:RFS}
\begin{algorithmic}[1]
\Require $y \in \mathbb{R}^K$, $X_0 \in \mathbb{R}^{K \times m}$, $\lambda \geq 0$, $y \neq 0$
\State $Scale_j \gets 1/\|X_0[:,j]\|_2$, \; $X \gets X_0 \cdot Scale$ \Comment{normalize columns}
\State $W \gets 0 \in \mathbb{R}^{m \times m}$, \; $Res \gets 0 \in \mathbb{R}^{m+1}$, \; $big \gets \{0,\dots,m-1\}$
\State $W[:,0] \gets \argmin_\psi (\|y - X\psi\|_2^2 + \lambda\|\psi\|_2^2)$, \; $Res[0] \gets \|y - XW[:,0]\|_2^2 + \lambda\|W[:,0]\|_2^2$
\For{$j = 0, \dots, m-2$}
    \State $k_{\min} \gets \argmin_{k \in big} |W[k,j]|$, \; $big \gets big \setminus \{k_{\min}\}$
    \State $W[big,j+1] \gets \argmin_\psi \|y - X[:,big]\,\psi\|_2^2 + \lambda\|\psi\|_2^2$
    \State $Res[j+1] \gets \|y - XW[:,j+1]\|_2^2 + \lambda\|W[:,j+1]\|_2^2$
\EndFor
\State $Res[m] \gets \|y\|_2^2$
\State $j_{\max} \gets \argmax_j Res[j+1]/Res[j]$
\State \Return $W[:,j_{\max}] \cdot Scale$
\end{algorithmic}
\end{algorithm}

\subsection{Weak formalism}
In order for our PDE discovery methods to work better for noisy data, we introduce the weak form of our PDE, see also \cite{messenger2021weakODE,messenger2021weak}.
Let $\Psi: \mathbb{R}^d \to \mathbb{R}$ be a test function, meaning a smooth function of compact support. In our paper, we will use the positive part of a polynomial function extended by zero. We assume that we can compute derivatives of $\Psi$ analytically. The (spatial) convolution $\Psi * f$ of functions $f,\Psi:\mathbb{R}^d \to \mathbb{R}$ is given by
$$(\Psi * f)(x) = \int_{\mathbb{R}^d} f(\tilde x) \Psi(x - \tilde x) d\tilde x$$
Using integration by parts, we obtain $\Psi * (\partial_{x^i} f) = (\partial_{x^i} \Psi)* f$. It follows that for any linear differential operator $L$ with constant coefficients we have $\Psi * (Lf) = (L\Psi) * f$. In general, there are no such universal integration-by-parts formulas for nonlinear operators.

Assume now that we have noisy spatiotemporal data $\mathbf{U} = \{U(t_k,x_k)\} =\{ u(t_k, x_k)+ \ve(t_k, x_k)\}$ of a solution of $u_t = N(u) = \sum_{j}^m \xi_j N_j(u)$. We furthermore assume that this data come from a grid $(t_k, x_k) \in \mathbf{t} \times \mathbf{x} \subseteq [0,T] \times \mathbb{R}^d$, with respective mesh sizes $\Delta t$ and $\Delta x$. For spatiotemporal data $\mathbf{U} = \{U(t_k,x_k) \}$, we write
$$\Psi * \mathbf{U} = \left\{ (\Delta x)^d\sum_{\tilde x \in \mathbf{x}} \Psi(x_k-\tilde x) U(t_k, \tilde x) \right\},$$
for the corresponding discrete convolution. Recall that we represent $u_t = N(u)$ numerically as $\mathbf{U}_t = \mathbf{N} \xi$. By instead considering $\Psi * \mathbf{U}_t = (\Psi * \mathbf{N}) \xi$, representing $\Psi * u_t = \Psi * N(u)$, the averaging over the spatial dimensions seeks to minimize the effect of spatial noise. Furthermore, if we can write $N_j = L Q$, where $L$ is a linear differential operator and $Q$ has lower order than $N_j$, we can replace $\Psi * \mathbf{N}_j$ with $L\Psi * \mathbf{Q}$, where $L\Psi$ is computed analytically and $\mathbf{Q}$ computes $Q(u)$ numerically. Decreasing the order of derivatives for numerical derivatives further decreases the impact of noisy data.

\subsection{Libraries of differential operators}
All of our libraries will consist of polynomial differential operators of at most total order 4 and polynomial degree 3. All libraries we consider will be translation-equivariant, that is, we assume that any operator $N_j$ in our library satisfies $N_j(f) \circ T = N_j(f \circ T)$ for any smooth function $f: \mathbb{R}^d \to \mathbb{R}$ and translation $T(x) = x+b$. The latter property is equivalent to assuming that $N_j$ has no explicit dependence on the spatial variable $x \in \mathbb{R}^d$. A library which spans all such operators is the \emph{the polynomial jet}-library $\LibPj$ which consists of all operators on the form
$$f \mapsto \prod_{l=1}^n \partial_{x^1}^{\alpha_{l,1}}  \cdots \partial_{x^d}^{\alpha_{l,d}} f, \qquad 0 \leq n \leq 3, \quad \alpha_{l,s} \geq 0, \quad \sum_{l=1}^n \sum_{s=1}^d \alpha_{l,s} \leq 4.$$
We will also use this library for reference when comparing other libraries, as it is the most comprehensive one. We remark that the choice of total order 4 and polynomial degree 3 as upper bounds are done here in this experiment is made to set an upper limit to the size of libraries. Also, with higher total orders, there will also be dimension-dependent non-trivial linear dependence relations for equivariant operators, see \cite{ballerin2026equivariant}.

\subsubsection{Equivariant multigraph-based library}
The objective of the paper is to take advantage of the fact that we know that the operator $N$ is equivariant under all rigid motions. The first library $\LibM$ consists of all 20 operators in Table~\ref{tab:multigraph_operators}, which we will call \emph{the multigraph library}, since the basis directly corresponds to isomorphism classes of multigraphs. We will also use \emph{the modified multigraph} library $\LibMM$ where we use a different basis that allows us to take maximal advantage of integrations in the weak formulation. A list of members in $\LibMM$ is found in Table~\ref{tab:multigraph_operators2}. Notice that all the members of $\LibMM$ are on the form $LQ$ where $L$ is a linear operator, allowing for integration by parts, and $Q$ is an operator of total order at most 2. The construction of $\LibMM$ and a comparison between the two equivariant libraries are found in Section~\ref{sec:MandMMcorrespondence}.

\begin{table}[htbp]
  \centering
  \renewcommand{\arraystretch}{1.5} 
  \small
  \begin{tabular}{|rr|cccc|}
  \hline
    $p$ &  & $n=0$ & $n=1$ & $n=2$ & $n=3$          \\  \hline
    $0$ &  & $1$ & $f$ & $f^2$ & $f^3$  \\ \hline
    $1$ & IBP & & $\Delta f$ & $\Delta f^2$ & $ \Delta f^3$ \\
    &  & & & $\|\nabla f\|^2$ & $f \| \nabla f \|^2$  \\\hline
    $2$ & IBP & & $\Delta^2 f$ & $\Delta^2 f^2$ & $ \Delta^2 f^3$ \\
    & IBP & & & $\Delta \| \nabla f\|^2$ & $ \Delta ( f \|\nabla f\|^2)$ \\
    &  & & & $\| \nabla^2 f \|^2$ & $f\| \nabla^2 f\|^2$ \\
    &  & & & $(\Delta f)^2$ & $f (\Delta f)^2$ \\
    &  & & & &  $\nabla^2 f(\nabla f, \nabla f)$   \\
    &  & & & & $\| \nabla f \|^2  \Delta f$ \\
    \hline
  \end{tabular}
  \caption{A different basis of 20 $\Ort(d)$-equivariant polynomial operators better for integration by parts. This library is denoted by $\LibMM$. We mark with IBP whenever we can use integration by parts.}
  \label{tab:multigraph_operators2}
\end{table}

\subsubsection{Reference libraries}
While including more terms in a library makes it more expressive, it also opens up for more opportunities for false positives in the regression process. Therefore, in addition to the maximal PolyJet library $\LibPj$, we will also compare with other smaller libraries used earlier in the literature. The library $\LibFIND$ used in PDE-FIND \cite{rudy2017data} consists of operators on the form $f \mapsto f^n Lf$ where $L= \partial_{x^1}^{\alpha_1} \cdots \partial_{x^d}^{\alpha_d}$ is a linear differential operator, while the library $\LibWS$ used for weak SINDy \cite{messenger2021weak} considers only operators on the form $L f^n$ where $L$ is a linear differential operator. Notice that the latter is chosen so that one can always remove all numerically computed derivatives in the weak formulation using integration by parts. We restrict ourselves to the terms having the correct total order and polynomial degree. Due to the restrictions of these libraries, they will not be able to represent all of the PDEs that we consider in our experiments.

\begin{table}[htbp]
  \centering
  \small
  \begin{tabular}{l|cccccl}
    Equation & $\LibFIND$ & $\LibWS$ & $\LibPj$ & $\LibM$ & $\LibMM$ \\ \hline
    Allen--Cahn \eqref{AC}                & \checkmark  & \checkmark  & \checkmark   & \checkmark  & \checkmark  \\
    Swift--Hohenberg \eqref{SH}           & \checkmark  & \checkmark  & \checkmark   & \checkmark  & \checkmark  \\
    Cahn--Hilliard \eqref{CH}             & $\times$  & \checkmark  & \checkmark   & \checkmark  & \checkmark  \\
    Porous medium \eqref{PM}              & $\times$  & \checkmark  & \checkmark   & \checkmark  & \checkmark  \\
    Kuramoto--Sivashinsky \eqref{KS}      & $\times$  & $\times$ & \checkmark   & \checkmark  & \checkmark  \\
  \end{tabular}
  \caption{An overview of which equations the different libraries can represent. Unlike the 1-dimensional case, in which we can define $v = \partial_x u$, then $v_t = - v v_x - \Delta v - \Delta^2 v$ so that the KS equation is representable in $\LibFIND$ as seen in \cite{rudy2017data}, the 2D case is spanned by neither $\LibFIND$ nor $\LibWS$. Thus, $\LibPj$ is the only non-equivariant library under consideration that can represent it. }
  \label{tab:library_span}
\end{table}

\subsection{Derivatives and noise} \label{sec:Derivatives}
Building the library of differential terms $\mathbf{N}$ requires numerical differentiation of the observed data. There are several choices in the literature which each come with their own strengths and weaknesses, especially when noise is introduced into the mix. Following the approach discussed by \cite{brunton2016discovering}, throughout our work we choose and compare two complementary approaches to computing derivatives, which we call ``backend'': \emph{spectral} and \emph{polynomial} backends. In all cases, computing differential operators is performed in terms of simpler derivatives, and for all libraries each differential term is written in terms of these simpler derivatives.

\subsubsection{Spectral backend. }
Let $\omega$ denote the variable in Fourier space and let $\hat f = \mathcal{F}(f)$ denote the Fourier transform of a function $f:\mathbb{R}^d \to \mathbb{R}$. Differentiation can be performed in Fourier space via the property
\[
\partial_{x^1}^{\alpha_1} \cdots \partial_{x^d}^{\alpha_d} f \;\longleftrightarrow\; (i \omega^1)^{\alpha_1} \cdots (i \omega^d)^{\alpha_d} \, \widehat{f}(\omega).
\]
A single forward FFT, a multiplication by the appropriate polynomial in $i \omega^1$, $\dots,$ $i \omega^d$, and one inverse FFT yield the derivative at every grid point simultaneously. The method is extremely fast and machine-precision exact for band-limited, periodic, noiseless data. However, spectral differentiation is highly susceptible to observational noise and is therefore not a good choice when it comes to discovering equations in noisy regimes, which is unfortunately the default setting for real-life applications.

\subsubsection{Polynomial backend. }
An alternative to spectral derivatives is fitting a polynomial stencil on the data and differentiating the polynomial instead. In particular, given a parameter $r\in\mathbb N_{>0}$ we can fit a polynomial of a given total degree to a $(2r+1)^d$ periodic neighborhood of each grid point $x_k$ by least squares as
\[f(x) \approx \sum_{\alpha_1+ \cdots+ \alpha_d \leq \text{deg}} c_{\alpha_1, \dots, \alpha_d}\, (x^1-x_k^1)^{\alpha_1} \cdots  (x^d - x_k^d)^{\alpha_d}, \qquad \text{close to $x_k$.}\]
The required partial derivative $\partial_{x^1}^{\alpha_1} \cdots \partial_{x^d}^{\alpha_d} f(x_k)$ is extracted analytically from the fitted coefficient $c_{\alpha_1, \dots, \alpha_d}$ via the formula $\alpha_1! \cdots \alpha_d!\,c_{\alpha_1, \dots,\alpha_d}$. Fitting a low-degree polynomial to many neighboring points acts as implicit smoothing before differentiation, trading a small approximation bias for improved noise robustness at the expense of some precision. The default hyperparameters $(deg, r)$ which govern the degree of the polynomial and the stencil size are empirically chosen. A few examples can be consulted in Table~\ref{tab:polynomial_derivatives_hyperparam}.

\begin{table}[]
    \centering
    \begin{tabular}{c|c|c} 
         Term & Order & ($deg$, $r$) \\ \midrule
         Gradient & 1 & (9,5) \\ 
         Hessian & 2 & (9,5) \\
         Laplacian & 2 & (9,5) \\
         3rd order partials & 3 & (9,5) \\
         Bilaplacian & 4 & (11,6) \\
         Trilaplacian & 6 & (13,7) \\ 
    \end{tabular}
    \caption{Order of the polynomial and size of grid for different differential terms, as implemented in our code for dimension $2$.}
    \label{tab:polynomial_derivatives_hyperparam}
\end{table}

\begin{remark}
Other differential backends can be implemented and tested, however this is out of the scope of this work. We refer to the supplement of \cite{rudy2017data} for a complete discussion.
\end{remark}

\subsubsection{Type of noise used} \label{sec:Noise}

We consider additive, independent, identically distributed Gaussian observation noise. After a clean simulation produces snapshot array $\mathbf{u} \in \mathbb{R}^{T \times N \times N}$, noisy snapshots are generated as
\[
\mathbf{U}(t_k,x_k) = \mathbf{u}(t_k,x_k) + \sigma \cdot \eta(t_k,x_k), \qquad \eta(t_k,x_k) \sim \mathcal{N}(0,1) \;\text{i.i.d.},
\]
where $\sigma = \nu \cdot \operatorname{std}(\mathbf{u})$ and $\nu \in \{0.0, 0.02, 0.05\}$ is the noise fraction. The noise is added independently to each saved snapshot, corresponding to uncorrelated measurement error at each observation time. This is the standard noise model used in the PDE-FIND and WSINDy literature, facilitating direct comparison. In comparison with other types of noise, we will call this noise of type $(G)$.

In addition to isotropic noise, we also experiment with \emph{spatially-varying} noise profiles to test whether the equivariant libraries provide a robustness advantage when the noise itself breaks symmetry under rigid motion. Four profile types are considered:
\begin{enumerate}[$\bullet$]
\item \textbf{Radial decay $(R_{\downarrow})$:} $P(r) \propto \exp(-r^2/(2s^2))$ - with $r = \| x\|_2$, the noise is maximal at the origin and decaying outward (rotational symmetry with respect to the origin).
\item \textbf{Radial growth $(R_{\uparrow})$:} $P(r) \propto 1 + A(1 - \exp(-r^2/(2s^2)))$. Again with $r = \|x\|_2$, the noise is now minimal at the origin, growing outward (rotational symmetry with respect to the origin).
\item \textbf{Axial decay/growth, $(A_{\downarrow})$/$(A_{\uparrow})$:} we use the same profiles as for radial decay/growth, but now with $r = |x^1|$, meaning that we only have changes in noise levels along the directions of one axis and without any rotational symmetry.
\end{enumerate}
All profiles are normalized so that $\mathbb{E}[\|P\|_F^2] = 1$, preserving the same total noise power as the isotropic case for a given noise fraction $\nu$. 

\subsection{Evaluation of the result}
Let $w$ be a solution using one of the methods mentioned above, relative to a library $\scrL =\{N_1, \dots, N_m\}$. Let $\xi$ be the true solution. We want to evaluate how far our estimate is from the true result.

The central evaluation of a result is the number of false negatives $FN = \# \{ j \, :\, w_j = 0, \xi_j \neq 0\}$ and false positives $FP = \# \{ j \, :\, w_j \neq 0, \xi_j = 0\}$. False negatives in particular show that the method has not able to identify the terms of the equations. We note that the number of false positives and false negatives depend on our library $\scrL$ of differential operators. We say that the equation is \emph{recovered} if there are neither false positives, nor false negatives. When we can then also look at the errors in the coefficients
$$e_{coef} = \left\{ \left| \frac{w_j-\xi_j}{\xi_j} \right| \,: \, \xi_j \neq 0\right\}.$$
Evaluating the mean and standard deviation of these values was used in PDE-FIND \cite{rudy2017data} and is reported along with the standard deviation of the differences. The advantage of this error is that it does not depend on the data. However, the error depends on the library and it also consider all coefficients equal, regardless of their impact on the solution.

We will introduce evaluations that are related to the data itself. From our underlying dataset $\mathbf{U}$, we consider $(\mathbf{U}_t, \mathbf{N})$, computed according the library $\LibPj$. We can then measure the
$$E_{data}(w) = \frac{\| \mathbf{U}_{t} - \mathbf{N} w\|_2}{\|\mathbf{U}_t\|_2}.$$
This error is a measure of accuracy in cases of true PDE-discovery where we do not know the actual solution. However, such an error will use noisy data, which means that a solution that is less affected by noise can have a higher error.
In our experiment, we have access to the true solution $\xi$, as well as clean data $\mathbf{u}$. We consider the corresponding matrix $\mathbf{n}$ using the operators in $\LibPj$ computed with spectral derivatives. We can then use the error
$$E_{coef}(w) = \frac{\|\mathbf{n}(\xi-w)\|_2}{\|\mathbf{n}\xi\|_2}$$
We remark that $E_{data}(w) \leq E_{base} + \frac{\|\mathbf{n} \xi\|_2}{\|\mathbf{U}_t\|_2} \left( E_{noise}(w) +  E_{coef}(w) \right)$ with
$$E_{base} = \frac{\|\mathbf{U}_t - \mathbf{N} \xi\|_2}{\| \mathbf{U}_t\|_2} \, , \qquad  E_{noise}(w) = \frac{\|(\mathbf{N}-\mathbf{n})(\xi-w)\|_2}{\|\mathbf{n} \xi\|_2}.$$
Here, $E_{base}$ is the error that the true solution will have relative to the noisy data and can be considered as a measure of the difficulty of discovering the PDE, while $E_{noise}(w)$ can be considered as measuring how much the difference of $w$ from the true solution is because of the noisy data. We will use the PolyJet-library to compute $E_{data}$, $E_{base}$, $E_{noise}$ and $E_{coef}$ in all experiments, their values will be the same for all libraries up to numerical errors in computing the derivatives.

\section{Results}

The results presented herein are fully reproducible using the accompanying code, available at the public repository \texttt{equivariant-data-discovery} \footnote{https://github.com/ballerin/equivariant-data-driven-discovery}.

\subsection{Results in dimension two} \label{sec:ResultsTwoD}
We start by working with all the aforementioned equations in dimension two. We will consider the multigraph and modified multigraph libraries measured up against the libraries from PDE-FIND, weak SINDy and the PolyJet library. Since we are mainly preoccupied with the libraries themselves, we have tried several different backends when it comes to regression methods and ways of computing the derivatives, and focus on the ideal case for each library. We will go into more details on these choices in the ablation study in Section~\ref{sec:Ablation}. In general, it is always an advantage to use a weak formulation, though taking advantage of integration by parts is mainly beneficial at higher levels of noise. We remind the reader that while only the weak SINDy library $\LibWS$ can take advantage of integration by parts for all terms, all of the libraries can use integration by parts for some terms. Spectral derivatives are usually the choice that gives optimal results for clean data, and then derivatives obtained by local polynomial regressions generally perform better as the data gets more noisy.

An overview of the results of the experiments comparing the equivariant libraries to the non-equivariant libraries are found in Table~\ref{tab:Results2d}. In short, all libraries achieve recovery for most equations for $0 \%$, $2 \%$ and $5 \%$ noise with low errors. Recall that we define recovery as correctly identifying the terms of the equations (no false positives or negatives). The exceptions are all regimes of the Kuramoto--Sivashinsky (KS) equation, which are difficult as soon as there is noise present in the data. Recovery only differs in that only the equivariant libraries, both the multigraph $\LibM$ and the modified multigraph $\LibMM$ library, manage to recover the Cahn--Hilliard equation at $5 \%$-noise as well as the KS equation at $2 \%$ noise. In Table~\ref{tab:Outputs2d}, we display the concrete outputs of coefficient from the equivariant libraries and mean difference in coefficients.

The Allen--Cahn equation \eqref{AC} is the simplest benchmark in our suite, containing only pure monomials ($u$, $u^3$) and pure second derivatives ($u_{xx}$, $u_{yy}$).
Under the weak formulation, all method/algorithm pairs achieve exact support recovery  across all noise levels with low coefficient errors. We can consider the Allen--Cahn equation too simple to serve as a differentiating benchmark as all libraries and backends achieve recovery. The Swift--Hohenberg equation \eqref{SH} is more challenging than Allen--Cahn due to the presence of the fourth-order biharmonic operator $\Delta^2 u = u_{xxxx} + 2u_{xxyy} + u_{yyyy}$. All libraries still achieve recovery, but some only with the correct backend. The porous media equation \eqref{PM} is the best showcase of the intrinsic advantage of a sparse library. In particular, the modified multigraph library is able to represent the whole dynamics with only one term giving such a basis a disproportionate advantage. Both the standard multigraph basis and the modified one consistently outperform most non-equivariant methods, especially at higher noise regimes. The Cahn--Hilliard equation can be considered as a more challenging version of the Swift--Hohenberg equation, but since one can take advantage of integration by parts for all terms, it favors Modified-Multigraphs and WSINDy. That being said, only the equivariant libraries manage recovery at $5\%$. Finally, the Kuramoto--Sivashinsky (KS) the hardest benchmark, in which we only get recovery in the case of clean data. Here we are also only left with the PolyJet $\LibPj$ for non-equivariant libraries. The equivariant methods generally obtain lower errors and are the only ones that obtain recovery in the cellular regime with 2 \% noise.

\begin{table}[ht]
\centering
\small
\caption{The results below show the best equivariant result from the equivariant libraries measured against the best results of the non-equivariant libraries. The columns list if there is recovery for each method. In addition, the error $E_{coef}$ is reported for the equivariant method. The final column reports how much the best equivariant error is below or above non-equivariant method. Negative percentage means that the equivariant method has the lowest error.}
\label{tab:Results2d}
\begin{tabular}{lc|lcc|lc||c}
\toprule
Equation & Noise & Best equiv. & Recovery & $E_{coef}$ & Best non-equiv. & Recovery & $
\% \Delta E_{coef}$ \\
\midrule
Allen--Cahn & $0 \%$ & Mod.-Multigr. & $\checkmark$ & $0.0048$ & PDE-FIND & $\checkmark$ & $-0.4 \%$  \\
 & $2 \%$ & Mod.-Multigr. & $\checkmark$ & $0.0048$ & WSINDy & $\checkmark$ & $-0.2 \%$  \\
 & $5 \%$ & Mod.-Multigr. & $\checkmark$ & $0.0063$ & PolyJet & $\checkmark$ & $+0.3 \%$  \\
\midrule
Swift--Hohenberg & $0 \%$ & Multigraph & $\checkmark$ & $0.0224$ & PDE-FIND & $\checkmark$ & $-2.8 \%$  \\
 & $2 \%$ & Multigraph & $\checkmark$ & $0.0278$ & PolyJet & $\checkmark$ & $+17.5 \%$  \\
 & $5 \%$ & Multigraph & $\checkmark$ & $0.0885$ & WSINDy & $\checkmark$ & $-0.7 \%$  \\
\midrule
Porous medium equation & $0 \%$ & Mod.-Multigr. & $\checkmark$ & $0.0007$ & WSINDy & $\checkmark$ & $-5.8 \%$  \\
& $2 \%$ & Mod.-Multigr. & $\checkmark$ & $0.0000$ & WSINDy & $\checkmark$ & $-98.7 \%$  \\
 & $5 \%$ & Mod.-Multigr. & $\checkmark$ & $0.0013$ & WSINDy & $\checkmark$ & $-77.1 \%$  \\
\midrule
Cahn--Hilliard & $0 \%$ & Mod.-Multigr. & $\checkmark$ & $0.0136$ & WSINDy & $\checkmark$ & $-15.7 \%$  \\
 & $2 \%$ & Mod.-Multigr. & $\checkmark$ & $0.0361$ & WSINDy & $\checkmark$ & $+4.9 \%$  \\
 & $5 \%$ & Mod.-Multigr. & $\checkmark$ & $0.1417$ & WSINDy & $\times$ &   \\
\midrule
KS (cellular) & $0 \%$ & Multigraph & $\checkmark$ & $0.0027$ & PolyJet & $\checkmark$ & $-12.2 \%$  \\
 & $2 \%$ & Multigraph & $\checkmark$ & $0.2711$ & PolyJet & $\times$ &   \\
 & $5 \%$ & Mod.-Multigr. & $\times$ &  & PolyJet & $\times$ &   \\
\midrule
KS (weakly chaotic) & $0 \%$ & Mod.-Multigr. & $\checkmark$ & $0.0021$ & PolyJet & $\checkmark$ & $-6.5 \%$  \\
 & $2 \%$ & Mod.-Multigr. & $\times$ &  & PolyJet & $\times$ &  \\
\midrule
KS (strongly chaotic) & $0 \%$ & Multigraph & $\checkmark$ & $0.0017$ & PolyJet & $\checkmark$ &  $-12.7 \%$ \\
 & $2 \%$ & Mod.-Multigr. & $\times$ &  & PolyJet & $\times$ &   \\
\bottomrule
\end{tabular}
\end{table}

\begin{table}[ht]
\centering
\small
\caption{Here we can see the true equations compared to the discovered equations found by the equivariant libraries at $2 \%$ noise, in order to give some context to the term $E_{coef}$. We also list the mean of coefficient error $e_{coef}$ along with the standard deviation. For the KS equation in the weakly chaotic regime, we fail to recover the equations, as we have a false positive constant term.}
\label{tab:Outputs2d}
\begin{tabular}{ll|c||c|c}
\toprule
Equation & & $N(f)$ & $E_{coef}$ & $e_{coef}$ \\
\midrule
AC & True & $0.0225\,\Delta f + f - f^3$ & 0.0048 & $0.57 \pm 0.10 \%$ \\
& Disc. & $(+0.02238)\,\Delta f  +(+0.99519)f + (-0.99318)f^3 $\\ \midrule
SH & True & $-0.7f - 2\,\Delta f - \Delta^2 f - f^3$ & 0.0278 & $0.87 \pm	0.75 \%$ \\
& Disc. & $(-0.70074)f + (-2.01422)\Delta f + (-1.00788)\Delta^2f+ (-1.01895)f^3$\\ \midrule
PME & True &  $f^2\Delta f + 2f\lVert\nabla f\rVert^2 = \frac{1}{3}\Delta f^3$ & 0.0000 & $0.002 \%$ \\ 
& Disc. & $(+0.33334)\Delta f^3$ \\ \midrule
CH & True &  $-\Delta f + \Delta f^3 - 0.04 \Delta^2 f$ & 0.0361 & $1.14 \pm	0.07\%$ \\
& Disc. &  $(-0.98808)\Delta f + (+0.98946)\Delta f^3 + (-0.03954)\Delta^2 f$\\ \midrule
KS & True & $-\frac{1}{2}\,\|\nabla f\|^2 - \Delta f - \Delta^2 f$ & 0.2711 & $38 \pm	16 \% $ \\
(cell.)& Disc. & $(-0.40074)\lVert\nabla f\rVert^2 + (-0.55605)\Delta f + (-0.50084)\Delta^2 f$\\ \midrule
KS & True & $-\frac{1}{2}\,\|\nabla f\|^2 - \Delta f - \Delta^2 f$  & 0.6262\\ 
(w. ch.)& Disc. & $(-0.36236)\lVert\nabla f\rVert^2 + (-0.68916)\Delta f  + (-0.66853)\Delta^2 f + (+1.01145)1$\\
\bottomrule
\end{tabular}
\end{table}

\subsection{Results with space-varying noise} \label{sec:SpaceVarying}
In our experiments in Section~\ref{sec:ResultsTwoD}, we are looking at Gaussian noise added at every point in spacetime to our data, labeled $(G)$. This means that we are considering an equation invariant under rigid motion of space along with a noise profile which is invariant under those symmetries as well. We also tested radial increasing/decreasing noise and axial increasing/decreasing noise as discussed in Section~\ref{sec:Noise}, to check how the different choices of libraries is affected by nonequivariant noise.

Recall that we are considering noise that is increasing/decreasing in radial directions $(R_{\uparrow}, R_{\downarrow})$ as well as axial increasing/decreasing in the $x$-direction $(A_{\uparrow}, A_{\downarrow})$. The cases where we have recovery are summarized in Table~\ref{tab:ResultsNoiseTypes}. The most challenging noise profile seems to be the radial decreasing noise. Under radial decreasing noise, none of the non-equivariant libraries of PDE-FIND, WSINDy or PolyJet recover the Swift--Hohenberg equation at $5 \%$, while it is recovered by the equivariant libraries. On the other hand, the usual multigraph library $\LibM$ does not recover the Cahn--Hilliard equation under 5 \% radial decreasing noise, while all types of noise and all cases for the modified-multigraph library $\LibMM$ results in recovery. None of the non-equivariant libraries have recovery for any of the noise profiles for Cahn--Hilliard at $5 \%$. We have very few cases for recovery for the KS equation, but at $2 \%$ we have recovery for the equivariant libraries in one case in addition to the Gaussian noise.

In general, the higher order derivatives, the more vulnerable we are to noise. Hence, it is not surprising that weak SINDy is most robust of the non-equivariant methods, as it can avoid computing derivatives through integration by parts. On the other hand, PDE-FIND and PolyJet libraries compute derivatives in axial directions, so it is then understandable that they will be affected by the axial noise that is anisotropic. The modified multigraph library performs the best out of all the libraries, as it is able to take advantage of integration by parts for many of its terms, while not having derivatives going purely in the direction of one axis.

\begin{table}[ht]
\centering
\caption{The table shows the case when we manage to have recovery of the different libraries with respect to different types of noise; axial increasing ($A_{\uparrow}$), axial increasing ($A_{\downarrow}$), radially increasing $(R_{\uparrow})$, radially decreasing $(R_{\downarrow})$ and the Gaussian noise ($G$) of experiments in Section~\ref{sec:ResultsTwoD}. Cases where we have recovery are listed. Results are not reported for a library when it is not able to express the equation in question.}
\label{tab:ResultsNoiseTypes}
\begin{tabular}{lc|c|c||c|c|c|}
\toprule
Equation & Noise &  $\LibM$ & $\LibMM$ & $\LibFIND$ & $\LibWS$ & $\LibPj$ \\
\midrule
Allen--Cahn & $2 \%$ & All & All & All & All & All \\
 & $5 \%$ & All & All & All & All & All  \\
Swift--Hohenberg & $2 \%$ & All & All & All & All & All \\
& $5 \%$ & All & All & $G, A_{\uparrow}, R_{\uparrow}$  & $G, A_{\uparrow}, A_{\downarrow}, R_{\uparrow}$ & $G, A_{\uparrow}, R_{\uparrow}$ \\
Porous medium equation & $2 \%$  & All & All & & All & All \\
& $5 \%$  & All & All & & All & $G, R_{\uparrow}$ \\
Cahn--Hilliard & $2 \%$  & All & All & & All& All \\
& $5 \%$  & $G, A_{\uparrow}, A_{\downarrow}, R_{\uparrow}$ & All & & None & None \\
KS cellular & $2 \%$ & $G, A_{\uparrow}$ & $G, R_{\uparrow}$ & &  & None\\
& $5 \%$ &  None & None & & & None\\
KS w. chaotic & $2 \%$ & None & None & & & None \\
 & $5 \%$ & None & None & & & None \\
\bottomrule
\end{tabular}
\end{table}

\subsection{Result with increased dimension} \label{sec:IncreasedDimension}
We next wanted to look at the effect on increasing the spatial dimension on the results. Because of the increased computation time in for increased dimension, we will focus on only one backend. We will only use STRidge for regression, applied to the libraries $\LibM$, $\LibMM$, $\LibFIND$ and $\LibWS$. The number of candidate terms is listed in Table~\ref{tab:DimLibraries}. With the increased dimension, it is also more resource intensive to complete local polynomial regression for every derivative, so we are using only spectral derivatives. To compensate for this, we use a higher discretization level for the tolerance parameter $tol >0$ and the ridge parameter $\lambda >0$ for STRidge in order to have hope of still recovering the equations with higher noise levels at the cost of longer computation times.

\begin{table}[ht]
\centering
\caption{Number of candidate terms in the libraries we tested, depending on the dimension $d$.}
\label{tab:DimLibraries}
\begin{tabular}{l|cccc}
\toprule
dim $d$ & $\LibM$ & $\LibMM$ & $\LibFIND$ & $\LibWS$  \\
\midrule
2 & 20 & 20 & 60 & 46\\
3 & 20 & 20 & 140 & 106 \\
4 & 20 & 20 & 280 & 211\\
\bottomrule
\end{tabular}
\end{table}

We tested three equations: the Allen--Cahn, Swift--Hohenberg and the porous medium equation. Our expectation was that with the increased dimensions, errors would increase and recovery would become more difficult. In particular, we expected the performance of $\LibFIND$ and $\LibWS$ to be affected by the larger number of candidate terms that comes with increased dimension.  Surprisingly, our main observation is that the increase in dimension does not impact much the actual performance of the method. We have recovery for all libraries, for all equations and for all levels of noise except for $\LibM$ for the Allen-Cahn equation at $5\%$. Furthermore, errors are comparable for the dimensions 2, 3 and 4. See Table~\ref{tab:ResultsHidherD} for the errors listed for the different libraries according to dimensions. Further research is needed to determine how the libraries differ at more difficult higher-dimensional tasks.

\begin{table}[ht]
\centering
\small
\caption{Below we find the $E_{coef}$ reported for all experiments at $2 \%$. }
\label{tab:ResultsHidherD}
\begin{tabular}{lc|cccc}
\toprule
Equation &  Dim & $\LibM$ & $\LibMM$ & $\LibFIND$ & $\LibWS$  \\
\midrule
             & 2 & 0.0057 & 0.0057 & 0.0059 & 0.0059 \\
 Allen--Cahn & 3 & 0.0070 & 0.0070 & 0.0071 & 0.0071\\
             & 4 & 0.0065 & 0.0065 & 0.0064 & 0.0064\\
\midrule
                  & 2 & 0.0067 &	0.0067 & 0.0058 & 0.0058 \\
 Swift--Hohenberg & 3 & 0.0184 & 0.0184& 0.0227 & 0.0227 \\
                  & 4 & 0.0264 & 0.0264 & 0.0268 & 0.0268\\
\midrule
                       & 2 & 0.0042 & 0.0035& --- & 0.0089 \\
Porous medium equation & 3 & 0.0045 & 0.0006 & --- & 0.0012\\
                       & 4 & 0.0014 & 0.0016 & --- & 0.0017 \\
\bottomrule
\end{tabular}
\end{table}

Nevertheless, we do see that the multigraph and modified multigraph library require much less computational time and RAM, compared to non-equivariant libraries. See details in Table~\ref{tab:computational_cost} for the Allen-Cahn equations. Hence, choosing equivariant libraries could be necessary if computational resources are limited.

\begin{table}[ht]
\centering
\caption{Computational costs for different libraries, depending on the dimension $d$. Data reported is for Allen--Cahn equation, noiseless data and pointwise formulation with spectral derivatives. Grid of 32 points per dimension, and library matrix is subsampled to 10,000 rows. Subsampling makes the problem tractable by cutting the growth in library size with dimension. Computations in dimension 5 are performed only for the multigraph libraries, as they are the the only ones that required more than the maximum available RAM which was allocated to 96GB.}

\label{tab:computational_cost}
\begin{tabular}{ll|ccc}
\toprule
Library & $d$ & Build time $\mathbf{N}$ (s) & STRidge (s) & RAM (GB) \\
\midrule
\multirow{3}{*}{$\LibFIND$}
 & 2 & 0.03 & 0.005 & 0.06 \\
 & 3 & 0.84 & 0.013 & 0.67 \\
 & 4 & 35.7 & 0.027 & 15.9 \\
\midrule
\multirow{3}{*}{$\LibWS$}
 & 2 & 0.10 & 0.003 & 0.06 \\
 & 3 & 1.66 & 0.012 & 0.67 \\
 & 4 & 61.4 & 0.039 & 15.9 \\
\midrule
\multirow{3}{*}{$\LibPj$}
 & 2 & 0.09 & 0.013 & 0.07 \\
 & 3 & 1.85 & 0.036 & 1.40 \\
 & 4 & 75.5 & 0.507 & 38.8 \\
\midrule
\multirow{4}{*}{$\LibM$}
 & 2 & 0.07 & 0.002 & 0.04 \\
 & 3 & 0.71 & 0.001 & 0.15 \\
 & 4 & 17.8 & 0.002 & 1.61 \\
 & 5 & 224  & 0.002 & 20.4 \\
\midrule
\multirow{4}{*}{$\LibMM$}
 & 2 & 0.10 & 0.001 & 0.04 \\
 & 3 & 0.94 & 0.001 & 0.15 \\
 & 4 & 22.6 & 0.002 & 1.64 \\
 & 5 & 304  & 0.002 & 20.6 \\
\bottomrule
\end{tabular}
\end{table}

\subsection{Ablation study} \label{sec:Ablation}

We conducted a broad ablation study across hyperparameters, libraries, and regression methods for the 2D case with space-constant Gaussian noise.  Each equation was evaluated with all relevant libraries at 0\%, 2\%, and 5\% observation noise, under every combination of backend (spectral, polynomial), formulation (pointwise, weak direct, weak with IBP), and regression method (STRidge, MSTLS, RFE).  We report results at the highest noise level for which at least one library--hyperparameter combination achieved exact support recovery: 5\% for Allen--Cahn, Swift--Hohenberg, Porous Medium, and Cahn--Hilliard; 2\% for Kuramoto--Sivashinsky in the cellular regime; and 0\% for Kuramoto--Sivashinsky in the weakly chaotic regime.  No configuration recovered the correct support for Kuramoto--Sivashinsky in the chaotic regime.

With the exception of Allen--Cahn, for which most libraries perform well, at least one configuration of hyperparameters for multigraph-based libraries consistently outperforms the alternatives across all other equations. The tables report for each equation all the best combinations of library--regressor, and the hyperparameter combination which achieved the best result. The tables are listed in Appendix~\ref{sec:TablesAblation}, with the complete tables for all hyperparameters combinations are available in the companion repository.

\section{Conclusion}

We have introduced libraries of candidate terms for PDE discovery that are equivariant under rigid motions, built on the classification of such operators in \cite{ballerin2026equivariant}, and tested them against established non-equivariant libraries on a suite of five benchmark equations exhibiting different behaviours. Our observation is that leveraging symmetry to impose equivariance at the level of the library can be a productive design principle for sparse-regression-based discovery.

Three observations stand out. First, on the equations where the discovery problem is genuinely difficult (Cahn--Hilliard, Kuramoto--Sivashinsky), the equivariant libraries recover the correct support in strictly more cases than their non-equivariant counterparts: the multigraph and modified multigraph libraries are the only ones to recover the Cahn--Hilliard equation at 5\% noise, and the only ones to recover the Kuramoto--Sivashinsky equation in the cellular regime at 2\% noise. Second, this advantage is amplified when the observational noise itself fails to respect the symmetry of the dynamics. Third, the number of equivariant candidate terms is constant in the spatial dimension which translates directly into lower computational cost as the dimension grows, potentially turning an intractable problem into a tractable one.

The experiments concerning the extension to higher dimensions produced a result we did not anticipate: we had expected accuracy to degrade with dimension for classical libraries, on the grounds that more candidate terms mean more opportunities for false positives. In practice, all libraries we tested retained comparable coefficient errors, and the equivariant advantage in this regime is therefore one of computational cost and tractability rather than of accuracy.

Our experiments were conducted on synthetic data from scalar equations on flat, periodic domains, with candidate terms restricted to total order at most four and polynomial degree at most three. Extending the approach in any of these directions is natural future work. The classification in \cite{ballerin2026equivariant} covers constant curvature spaces, which would permit discovery on spheres and hyperbolic domains where the equivariance assumption is arguably more informative than in the flat case.

\bibliographystyle{abbrv}
\bibliography{Bibliography}

\appendix
\section{Construction of the modified multigraph library}
\label{sec:MandMMcorrespondence}
For the weak formalism, we want to take advantage of the fact that we know the derivatives of the test functions explicitly to do integration by parts whenever possible. Hence, we want to consider which equivariant operators $P$ we can write as $P(f) = L Q(f)$ where $L$ is a linear differential operator and $Q$ is a non-linear operator, both equivariant under rotations and reflections. Such equivariant linear operators $L$ are always polynomials in the Laplacian by \cite{helgason1959differential}. Since we are limiting ourselves to total degree 4, we know that $L$ has to be a second order polynomial in the Laplacian. We consider each power of the Laplacian separately,

If $L = \Delta$, then $Q$ can have total order at most~2, and polynomial degree at most 3. The options for $Q$ are then $f$, $f^2$, $f^3$, $\Delta f$, $f \Delta f$, $f^2 \Delta f$, $|\nabla f|^2$ and $f |\nabla f|^2$.
We use the general equations that for a $C^2$-function $F$ of one of two variables:
\begin{align*}
\Delta F(u) &=F'(u) \Delta u  +  F''(u) \| \nabla u\|^2, \\
\Delta F(u,v) & = (\partial_u F)(u,v) \Delta u + (\partial_v F)(u,v) \Delta v \\
& \quad +(\partial_u^2 F)(u,v) \|\nabla u\|^2 + 2 (\partial_u \partial_v F)(u,v) \langle \nabla u, \nabla  v \rangle  + (\partial_v^2 F)(u,v) |\nabla v|^2.
\end{align*}
For the case when $L = \Delta^2$, then $Q$ has to have total degree 0. We use that for a $C^4$-function $F$, we have
\begin{align*}
\Delta^2 F(u) & = F^{(4)}(u) \|\nabla u\|^4 +  F'''(u)  ( 4(\nabla^2 u)( \nabla u, \nabla u)  +  2 \|\nabla u \|^2 \Delta u ) \\
&  \qquad + F''(u) (4 \langle \nabla \Delta u , \nabla u \rangle + 2  \|\nabla^2 u \|^2  +  (\Delta u )^2) + F'(u) \Delta^2 u.
\end{align*}
Applying these formulas, we obtain the following result.
\begin{lemma} \label{lemma:WeakForm}
All equivariant polynomial differential operators of total order $\leq 4$ and $
\leq 3$ are spanned by the terms in Table~\ref{tab:multigraph_operators2} for the library $\LibMM$. We can convert terms from $\LibMM$ to $\LibM$ through Table~\ref{tab:LMtoLMM}.
\end{lemma}

\begin{table}[htbp]
  \centering
  \small
  \begin{tabular}{r|cc|cccccc} 
  \toprule
     & $f\Delta f$ & $\| \nabla f \|^2$
     & $f \Delta^2 f$ & $\langle \nabla \Delta f, \nabla f \rangle$ & $\| \nabla^2 f\|^2$ & $(\Delta f)^2$ &&
     \\  \midrule
    $\Delta f^2$ &2&2 &&&& && \\
    $\Delta \| \nabla f\|^2$ && &&2&2& && \\
    $\Delta^2 f^2$ && &2&8&4&2 && \\ \toprule
    & $f^2 \Delta f$ & $f \| \nabla f\|^2$
     & $f^2 \Delta^2 f$ & $f\langle \nabla \Delta f, \nabla f \rangle$ & $f \| \nabla^2 f\|$ & $f(\Delta f)^2$ & $\nabla^2 f(\nabla f, \nabla f)$ & $\|\nabla^2 f\|^2 \Delta f$ \\ \midrule
    $\Delta f^3$ & 3 & 6 &&&& && \\
    $\Delta (f\| \nabla f\|^2)$ && &&2&2&&4&1 \\
    $\Delta^2 f^3$ &&& 3&24&12&6&24&12 \\
    \bottomrule
  \end{tabular}
  \caption{The table shows the change of basis for elements in $\LibMM$ to $\LibM$.}
  \label{tab:LMtoLMM}
\end{table}

\section{Covariance between columns and conditioning}
 \label{sec:Covariance}

Ridge regression is employed in this work in order to improve the conditioning of the matrix representing the terms. In particular, experiments show that multicollinearity in the matrix $\mathbf{N}$ consisting of numerically computed nonlinear differential operators from the observed data $\mathbf{U}$ plagues recovery. We make this claim precise in this section, following \cite[Chapter~3.3]{james2013introduction}.

Pairwise correlations between columns can be detected in $\mathbf{N}$ by examining $\max_{i \neq j} |\text{cor}(\mathbf{N}_i,\mathbf{N}_j)|$. In order to detect collinearity involving multiple columns of $\mathbf{N}$, let $\widehat{\mathbf{N}}_{j}$ denote the matrix $\mathbf{N}$ with the column $\mathbf{N}_j$ removed. We can compute the \emph{variance inflation factor}

\[\mathrm{VIF}_j = \frac{1}{1-R^2_j}, \qquad R^2_j = 1 - \frac{\mathbf{N}_j (I - \widehat{\mathbf{N}}_{j} (\widehat{\mathbf{N}}_{j} \widehat{\mathbf{N}}_{j}^T)^{-1} \widehat{\mathbf{N}}_{j}^T ) \mathbf{N}_j}{\| \mathbf{N}_j-mean(\mathbf{N}_j)\|_2^2} \qquad j=1, \dots, m.
\]
In other words, $R_j^2$ is the $R^2$ of an OLS regression of $\mathbf{N}_j$ on the remaining columns of $\mathbf{N}$, and $\mathrm{VIF}_j$ measures the factor by which this correlation inflates the variance of the estimated regression coefficient of $\mathbf{N}_j$.

In addition to $\mathrm{VIF}$, we report the condition number of the column-normalized library matrix, $\kappa_{\mathrm{norm}} = \kappa(\tilde{N})$, where $\tilde{N}$ is obtained from $N$ by scaling each column to unit $\ell^2$-norm.

In Table~\ref{tab:collinearity_libraries}, we have computed these numbers for a selection of different equations and libraries, illustrating that multicollinearity is a problem that affects PDE recovery across all configurations. Across all configurations (equation, library, formulation, noise), no combination yields a design matrix that is well-posed with respect to multicollinearity. Moreover, the maximum absolute pairwise correlation never falls below 0.919, with a median of 0.990.

Surprisingly, library conditioning is not the predominant factor affecting recovery. Among the equations under consideration, Allen--Cahn and Cahn--Hilliard exhibit comparable values for VIF and condition number $\kappa_{\mathrm{norm}}$, but behave very differently during recovery in which the former appears as an easy benchmark for all libraries, whereas the latter is where the difficulty appears for standard non-equivariant libraries. Similar conclusions can be drawn comparing the Swift--Hohenberg and Porous Medium equations.

\begin{table}[ht]
\centering
\caption{Multicollinearity metrics across equations and libraries under zero-noise conditions. All differential operators are computed via the spectral backend.}
\label{tab:collinearity_libraries}
\begin{tabular}{llrrrr}
\toprule
Equation & Library & $\max|\mathrm{corr}|$ & max VIF & mean VIF & $\kappa_{\mathrm{norm}}$ \\
\midrule
  \multirow{5}{*}{AC}  & $\LibFIND$ & 0.919 & 15.3 & 5.2 & 13.4 \\
    & $\LibWS$ & 0.919 & 10.5 & 3.1 & 9.1 \\
    & $\LibM$ & 0.934 & 14.4 & 8.1 & 13.6 \\
    & $\LibMM$ & 0.945 & 377 & 48.3 & 57.9 \\
    & $\LibPj$ & 0.919 & 13.7 & 4.5 & 16.4 \\
\midrule
  \multirow{5}{*}{SH}  & $\LibFIND$ & 0.997 & 6,776 & 718 & 456 \\
    & $\LibWS$ & 0.996 & 2,108 & 177 & 151 \\
    & $\LibM$ & 0.999 & 460,643 & 46,117 & 2,562 \\
    & $\LibMM$ & 0.998 & 443,470 & 94,131 & 2,752 \\
    & $\LibPj$ & 0.998 & 551,901 & 10,873 & 4,283 \\
\midrule
  \multirow{4}{*}{PME}  & $\LibWS$ & 0.980 & 1,840 & 111 & 147 \\
    & $\LibM$ & 0.980 & 500,735 & 145,898 & 5,198 \\
    & $\LibMM$ & 0.980 & 512,072 & 156,192 & 4,914 \\
    & $\LibPj$ & 0.980 & 8,203,610 & 354,594 & 24,366 \\
\midrule
  \multirow{4}{*}{CH}  & $\LibWS$ & 0.975 & 206 & 11.4 & 37.6 \\
    & $\LibM$ & 0.975 & 456 & 59.2 & 74.5 \\
    & $\LibMM$ & 0.975 & 735 & 123 & 89.6 \\
    & $\LibPj$ & 0.975 & 461 & 19.0 & 86.9 \\
\bottomrule
\end{tabular}
\end{table}

Library multicollinearity therefore does not tell the whole story, and the actual terms in the libraries play a significant role in the task. Libraries computed with noisy data are generally affected by lower VIF. However, this is not a direct improvement but rather an artifact of the noise injecting independent variation into the columns.

Additionally, weak formulations are more poorly conditioned than their pointwise counterparts. One explanation is that smoothing against test functions removes the high-frequency content that distinguishes the library terms from one another. A second factor is that noise injects independent variation into each differentiated column, so that multicollinearity and condition number appear artificially low under noisy regimes.

\begin{table}[ht]
\centering
\caption{Percent change in $\mathrm{VIF}_{\max}$ when moving from the pointwise to the weak formulation, averaged over equations.  Positive values mean the weak formulation is more collinear.}
\label{tab:coll-pct-lib-max-vif}
\begin{tabular}{lrrrr}
\toprule
Library & \multicolumn{2}{c}{Clean} & \multicolumn{2}{c}{2\% noise} \\
\cmidrule(lr){2-3}\cmidrule(lr){4-5}
 & Mean & Median & Mean & Median \\
\midrule
$\LibFIND$ & +1\,197\% & +1\,197\% & +4\,488\% & +4\,488\% \\
$\LibWS$ & +893\% & +400\% & +7\,127\% & +3\,636\% \\
\quad +IBP & +892\% & +400\% & +7\,130\% & +3\,630\% \\
$\LibM$ & +2\,153\% & +572\% & +17\,287\% & +15\,669\% \\
$\LibMM$ & +1\,003\% & +215\% & +55\,229\% & +26\,347\% \\
\quad +IBP & +857\% & +17\% & +10\,474\% & +1\,741\% \\
$\LibPj$ & +23\,186\% & +824\% & +12\,766\% & +11\,984\% \\
\bottomrule
\end{tabular}
\end{table}

\begin{table}[ht]
\centering
\caption{Percent change in $\kappa_{\mathrm{norm}}$ when moving from the pointwise to the weak formulation, averaged over equations.  Positive values mean the weak formulation is more collinear.}
\label{tab:coll-pct-lib-cond-normalised}
\begin{tabular}{lrrrr}
\toprule
Library & \multicolumn{2}{c}{Clean} & \multicolumn{2}{c}{2\% noise} \\
\cmidrule(lr){2-3}\cmidrule(lr){4-5}
 & Mean & Median & Mean & Median \\
\midrule
$\LibFIND$ & +402\% & +402\% & +574\% & +574\% \\
$\LibWS$ & +334\% & +249\% & +779\% & +572\% \\
\quad +IBP & +334\% & +249\% & +779\% & +571\% \\
$\LibM$ & +365\% & +293\% & +2\,366\% & +2\,162\% \\
$\LibMM$ & +213\% & +134\% & +4\,157\% & +3\,327\% \\
\quad +IBP & +166\% & +34\% & +991\% & +461\% \\
$\LibPj$ & +324\,932\% & +352\% & +1\,929\% & +2\,004\% \\
\bottomrule
\end{tabular}
\end{table}

\section{Tables from the ablation study} \label{sec:TablesAblation}
Below we list the tables for all the equations with the performance of the different libraries and with different selections of regression methods. These are listed in Tables~\ref{tab:ac-5pct} to~\ref{tab:ks-wchaos-clean}. We report the result from the case of maximal noise for each equation, where we still have recovery in at least one case.

\begin{table}[H]
\centering
\caption{Allen-Cahn equation — 5\% noise}
\label{tab:ac-5pct}
\small
\begin{tabular}{llrrrr}
\toprule
Method & Best conf. & fp/fn & $E_\text{coeff}$ & Time (s) & vs.\ best \\
\midrule
PolyJet + MSTLS & poly., weak (no IBP), m=12 & 0/0 $\checkmark$ & 0.0063 & 0.9 & --- \\
PDE-FIND + MSTLS & poly., weak (no IBP), m=12 & 0/0 $\checkmark$ & 0.0063 & 0.5 & +0.0\% \\
Multigraphs + MSTLS & poly., weak (no IBP), m=12 & 0/0 $\checkmark$ & 0.0063 & 0.1 & +0.3\% \\
Mod. Multi. + MSTLS & poly., weak (no IBP), m=12 & 0/0 $\checkmark$ & 0.0063 & 0.1 & +0.3\% \\
Mod. Multi. + RFE & poly., weak (no IBP), m=8 & 0/0 $\checkmark$ & 0.0065 & 0.0 & +3.4\% \\
Multigraphs + RFE & poly., weak (no IBP), m=8 & 0/0 $\checkmark$ & 0.0065 & 0.0 & +3.4\% \\
PolyJet + STRidge & poly., weak (no IBP), m=12 & 0/0 $\checkmark$ & 0.0065 & 0.4 & +4.0\% \\
PDE-FIND + STRidge & poly., weak (no IBP), m=12 & 0/0 $\checkmark$ & 0.0065 & 0.3 & +4.0\% \\
WSINDy + STRidge & poly., weak (no IBP), m=12 & 0/0 $\checkmark$ & 0.0065 & 0.3 & +4.0\% \\
Mod. Multi. + STRidge & poly., weak (no IBP), m=12 & 0/0 $\checkmark$ & 0.0066 & 0.1 & +4.3\% \\
Multigraphs + STRidge & poly., weak (no IBP), m=12 & 0/0 $\checkmark$ & 0.0066 & 0.1 & +4.3\% \\
WSINDy + MSTLS & poly., weak (IBP), m=12 & 0/0 $\checkmark$ & 2.8244 & 0.4 & +44778.7\% \\
PDE-FIND + RFE & poly., weak (no IBP), m=4 & 0/2 $\times$ & 0.3299 & 0.1 & $\times$ \\
PolyJet + RFE & poly., weak (no IBP), m=4 & 0/2 $\times$ & 0.3299 & 0.1 & $\times$ \\
WSINDy + RFE & poly., weak (no IBP), m=4 & 0/2 $\times$ & 0.3299 & 0.1 & $\times$ \\
\bottomrule
\end{tabular}
\end{table}

\begin{table}[H]
\centering
\caption{Swift-Hohenberg equation — 5\% noise}
\label{tab:sh-5pct}
\small
\begin{tabular}{llrrrr}
\toprule
Method & Best conf. & fp/fn & $E_\text{coeff}$ & Time (s) & vs.\ best \\
\midrule
Multigraphs + STRidge & poly., weak (no IBP), m=12 & 0/0 $\checkmark$ & 0.0779 & 0.2 & --- \\
Mod. Multi. + STRidge & poly., weak (no IBP), m=12 & 0/0 $\checkmark$ & 0.0779 & 0.2 & +0.0\% \\
WSINDy + STRidge & poly., weak (IBP), m=12 & 0/0 $\checkmark$ & 0.0782 & 0.6 & +0.4\% \\
Multigraphs + RFE & poly., weak (no IBP), m=12 & 0/0 $\checkmark$ & 0.0799 & 0.1 & +2.5\% \\
Mod. Multi. + RFE & poly., weak (no IBP), m=12 & 0/0 $\checkmark$ & 0.0799 & 0.1 & +2.5\% \\
WSINDy + RFE & poly., weak (IBP), m=12 & 0/0 $\checkmark$ & 0.0803 & 0.1 & +3.1\% \\
PolyJet + RFE & poly., weak (IBP), m=12 & 0/0 $\checkmark$ & 0.0803 & 0.1 & +3.1\% \\
PDE-FIND + RFE & poly., weak (IBP), m=12 & 0/0 $\checkmark$ & 0.0803 & 0.1 & +3.1\% \\
PDE-FIND + STRidge & spectral, weak (IBP), m=12 & 0/0 $\checkmark$ & 0.0891 & 0.5 & +14.4\% \\
Mod. Multi. + MSTLS & spectral, weak (no IBP), m=12 & 2/0 $\times$ & 0.4611 & 0.3 & $\times$ \\
PDE-FIND + MSTLS & spectral, weak (IBP), m=12 & 3/0 $\times$ & 0.1693 & 1.2 & $\times$ \\
PolyJet + STRidge & spectral, weak (IBP), m=12 & 5/0 $\times$ & 0.1868 & 1.0 & $\times$ \\
Multigraphs + MSTLS & spectral, weak (no IBP), m=12 & 8/0 $\times$ & 0.4818 & 0.3 & $\times$ \\
PolyJet + MSTLS & spectral, weak (no IBP), m=12 & 21/0 $\times$ & 0.6253 & 2.7 & $\times$ \\
WSINDy + MSTLS & poly., weak (IBP), m=12 & 31/0 $\times$ & 86.7777 & 1.3 & $\times$ \\
\bottomrule
\end{tabular}
\end{table}

\begin{table}[H]
\centering
\caption{Porous Medium equation ($m = 3$) — 5\% noise}
\label{tab:pme-5pct}
\small
\begin{tabular}{llrrrr}
\toprule
Method & Best conf. & fp/fn & $E_\text{coeff}$ & Time (s) & vs.\ best \\
\midrule
Mod. Multi. + STRidge & spectral, weak (IBP), m=12 & 0/0 $\checkmark$ & 0.0013 & 0.1 & --- \\
Multigraphs + STRidge & poly., weak (no IBP), m=12 & 0/0 $\checkmark$ & 0.0020 & 0.1 & +55.2\% \\
Mod. Multi. + MSTLS & spectral, weak (IBP), m=12 & 0/0 $\checkmark$ & 0.0025 & 0.1 & +93.0\% \\
Mod. Multi. + RFE & spectral, weak (IBP), m=12 & 0/0 $\checkmark$ & 0.0025 & 0.0 & +93.0\% \\
Multigraphs + MSTLS & poly., weak (no IBP), m=12 & 0/0 $\checkmark$ & 0.0056 & 0.1 & +334.6\% \\
Multigraphs + RFE & poly., weak (no IBP), m=12 & 0/0 $\checkmark$ & 0.0056 & 0.0 & +334.6\% \\
WSINDy + STRidge & poly., weak (no IBP), m=12 & 0/0 $\checkmark$ & 0.0056 & 0.2 & +336.0\% \\
PolyJet + MSTLS & poly., weak (no IBP), m=8 & 0/0 $\checkmark$ & 0.0059 & 1.6 & +364.3\% \\
PolyJet + STRidge & poly., weak (no IBP), m=8 & 0/0 $\checkmark$ & 0.0193 & 0.6 & +1406.5\% \\
WSINDy + MSTLS & poly., weak (no IBP), m=4 & 0/0 $\checkmark$ & 0.0441 & 0.5 & +3346.2\% \\
WSINDy + RFE & poly., weak (no IBP), m=8 & 0/2 $\times$ & 0.5773 & 0.0 & $\times$ \\
PolyJet + RFE & poly., weak (no IBP), m=8 & 0/3 $\times$ & 0.6700 & 0.1 & $\times$ \\
\bottomrule
\end{tabular}
\end{table}

\begin{table}[H]
\centering
\caption{Cahn-Hilliard equation — 5\% noise}
\label{tab:ch-5pct}
\small
\begin{tabular}{llrrrr}
\toprule
Method & Best conf. & fp/fn & $E_\text{coeff}$ & Time (s) & vs.\ best \\
\midrule
Mod. Multi. + RFE & poly., weak (IBP), m=12 & 0/0 $\checkmark$ & 0.1271 & 0.7 & --- \\
Mod. Multi. + STRidge & poly., weak (IBP), m=12 & 0/0 $\checkmark$ & 0.1276 & 1.2 & +0.4\% \\
Mod. Multi. + MSTLS & poly., weak (IBP), m=12 & 0/0 $\checkmark$ & 0.1755 & 1.5 & +38.1\% \\
Multigraphs + RFE & poly., weak (IBP), m=8 & 0/0 $\checkmark$ & 0.5795 & 0.7 & +356.0\% \\
WSINDy + RFE & poly., weak (no IBP), m=12 & 2/0 $\times$ & 0.1894 & 0.9 & $\times$ \\
WSINDy + STRidge & poly., weak (IBP), m=12 & 2/0 $\times$ & 0.1916 & 2.0 & $\times$ \\
Multigraphs + STRidge & poly., weak (no IBP), m=8 & 7/0 $\times$ & 0.7128 & 1.2 & $\times$ \\
PolyJet + STRidge & poly., weak (IBP), m=8 & 8/0 $\times$ & 0.6992 & 4.6 & $\times$ \\
PolyJet + RFE & spectral, weak (no IBP), m=12 & 8/0 $\times$ & 1.5522 & 1.4 & $\times$ \\
WSINDy + MSTLS & spectral, weak (no IBP), m=4 & 8/0 $\times$ & 40.4156 & 3.2 & $\times$ \\
Multigraphs + MSTLS & poly., weak (IBP), m=8 & 9/0 $\times$ & 0.6860 & 1.5 & $\times$ \\
PolyJet + MSTLS & poly., weak (IBP), m=12 & 10/0 $\times$ & 0.5369 & 21.8 & $\times$ \\
\bottomrule
\end{tabular}
\end{table}

\begin{table}[H]
\centering
\caption{Kuramoto-Sivashinsky equation (cellular, $L=4\pi$) — 2\% noise}
\label{tab:ks-cell-2pct}
\small
\begin{tabular}{llrrrr}
\toprule
Method & Best conf. & fp/fn & $E_\text{coeff}$ & Time (s) & vs.\ best \\
\midrule
Multigraphs + STRidge & poly., weak (no IBP), m=12 & 0/0 $\checkmark$ & 0.2711 & 0.2 & --- \\
Multigraphs + MSTLS & poly., weak (no IBP), m=12 & 0/0 $\checkmark$ & 0.2712 & 0.2 & +0.1\% \\
Mod. Multi. + MSTLS & poly., weak (no IBP), m=12 & 0/0 $\checkmark$ & 0.2712 & 0.2 & +0.1\% \\
PolyJet + MSTLS & poly., weak (IBP), m=12 & 4/0 $\times$ & 0.3327 & 2.3 & $\times$ \\
PolyJet + STRidge & poly., weak (IBP), m=12 & 25/0 $\times$ & 0.3054 & 1.1 & $\times$ \\
Mod. Multi. + STRidge & poly., weak (no IBP), m=8 & 28/0 $\times$ & 0.4081 & 0.3 & $\times$ \\
Multigraphs + RFE & spectral, weak (no IBP), m=8 & 0/5 $\times$ & 0.7796 & 0.1 & $\times$ \\
Mod. Multi. + RFE & spectral, weak (no IBP), m=8 & 0/5 $\times$ & 0.7796 & 0.1 & $\times$ \\
\bottomrule
\end{tabular}
\end{table}

\begin{table}[H]
\centering
\caption{Kuramoto-Sivashinsky equation (weakly chaotic, $L=8\pi$) — 0\% noise}
\label{tab:ks-wchaos-clean}
\small
\begin{tabular}{llrrrr}
\toprule
Method & Best conf. & fp/fn & $E_\text{coeff}$ & Time (s) & vs.\ best \\
\midrule
Mod. Multi. + STRidge & spectral, weak (no IBP), m=12 & 0/0 $\checkmark$ & 0.0021 & 0.4 & --- \\
Multigraphs + STRidge & spectral, weak (no IBP), m=12 & 0/0 $\checkmark$ & 0.0021 & 0.3 & +0.0\% \\
Mod. Multi. + MSTLS & spectral, weak (no IBP), m=12 & 0/0 $\checkmark$ & 0.0022 & 0.5 & +3.6\% \\
Multigraphs + MSTLS & spectral, weak (no IBP), m=12 & 0/0 $\checkmark$ & 0.0022 & 0.3 & +3.6\% \\
Multigraphs + RFE & spectral, weak (no IBP), m=12 & 0/0 $\checkmark$ & 0.0022 & 0.1 & +3.6\% \\
Mod. Multi. + RFE & spectral, weak (no IBP), m=12 & 0/0 $\checkmark$ & 0.0022 & 0.2 & +3.6\% \\
PolyJet + STRidge & spectral, weak (no IBP), m=12 & 0/0 $\checkmark$ & 0.0023 & 1.4 & +6.9\% \\
PolyJet + MSTLS & spectral, weak (no IBP), m=12 & 0/0 $\checkmark$ & 0.0023 & 4.1 & +8.3\% \\
PolyJet + RFE & spectral, weak (no IBP), m=12 & 0/0 $\checkmark$ & 0.0023 & 0.2 & +8.3\% \\
\bottomrule
\end{tabular}
\end{table}

\end{document}